\documentclass[10pt,leqno]{amsart}
\usepackage{a4,latexsym,amssymb,amsfonts,amsmath, graphicx, float, floatflt}
\usepackage{enumerate}
\usepackage{url}
\usepackage{mathbbol}
\usepackage[latin1]{inputenc}
\numberwithin{equation}{section}
\newtheorem{theo}{Theorem}

\newtheorem{lemma}[theo]{Lemma}
\newtheorem{prop}[theo]{Proposition}
\newtheorem{cor}[theo]{Corollary}
\newtheorem{defi}[theo]{Definition}
\newtheorem{assum}[theo]{Assumption}

\theoremstyle{definition}

\newtheorem{rem}[theo]{Remark}
\newtheorem{rems}[theo]{Remarks}
\newtheorem{exa}[theo]{Example}

\numberwithin{theo}{section}
\def\me{\mathsf{e}}
\def\mv{\mathsf{v}}
\def\mE{\mathsf{E}}
\def\mG{\mathsf{G}}
\def\mV{\mathsf{V}}

\DeclareMathOperator{\Ker}{Ker}
\DeclareMathOperator{\Range}{Rg}
\DeclareMathOperator{\Id}{Id}

\DeclareMathOperator{\Real}{Re}

\def\K{\tilde{K}}

\def\I{\mathcal I}
\def\1{\mathbb 1}
\def\It{\tilde{\mathcal I}}
\def\C{\mathbb C}

\DeclareMathOperator{\diag}{diag}
\DeclareMathOperator{\Span}{span}

\def\@secappcntformat#1{%
 \ifappendix \rm\appendixname\ifoneappendix\else~\fi\fi 
 \ifoneappendix \else \csname the#1\endcsname\relax\fi
 \ifx\apphe@d\@empty \else .\fi\enskip   
}

\title[Symmetries in strongly coupled network equations]{Well-posedness and symmetries of\\ strongly coupled network equations}
\author{Stefano Cardanobile}
\author{Delio Mugnolo}
\author{Robin Nittka}
\address{Institut f\"ur Angewandte Analysis, Universit\"at Ulm, Helmholtzstra{\ss}e 18, D-89081 Ulm, Germany}
\email{stefano.cardanobile@uni-ulm.de}
\email{delio.mugnolo@uni-ulm.de}
\email{robin.nittka@uni-ulm.de}


\begin{document}
%
%
%
%
%

\begin{abstract}
We consider a class of evolution equations taking place on the edges of a finite network and allow for feedback effects between different, possibly non-adjacent edges. This generalizes the setting that is common in the literature, where the only considered interactions take place at the boundary, i.\,e., in the nodes of the network. We discuss well-posedness of the associated initial value problem as well as contractivity and positivity properties of its solutions. Finally, we discuss qualitative properties that can be formulated in terms of invariance of linear subspaces of the state space, i.\,e., of symmetries of the associated physical system. Applications to a neurobiological model as well as to a system of linear Schr\"odinger equations on a quantum graph are discussed.
\end{abstract}

\maketitle

\section{Introduction}

The mathematical analysis of elliptic operators acting on spaces of functions on networks was started by G.\ Lumer in~\cite{Lum80}--\cite{Lum80a}. It has been subsequently continued by many authors, both in mathematics (in the context of \emph{network diffusion problems}, see e.g.\ \cite{Rot84}--\cite{Bel85}--\cite{Nic87}) and in physics (leading to the theory of \emph{quantum graphs}, see e.g.\ \cite{ExnSeb89}--\cite{KotSmi97}--\cite{KosSch99}). 


A form of weak nonlocal interactions for evolutionary problems  over network-shaped structures has already been considered in e.g.~\cite{KosSch99}--\cite{Mug07a}. Additionally, we are interested in discussing systems of \emph{strongly} coupled evolution equations. Such couplings may correspond to the cases of either a phenomenological interaction among parts of the physical system (like in a certain neurophysical theory, which we briefly discuss in Section~5.1) or else as a form of external control (possibly with the aim of stabilization).

More precisely, we want to allow the evolution in a point of the network to depend nonlocally on those finitely many other points of the network $\mG$
that have same parametrization with respect to the network edges. In other words, we will discuss the strongly coupled elliptic operator defined by
\begin{equation}\label{parab}
(Au)_j(x):=\sum_{i=1}^m\frac{\partial}{\partial x}\left(c_{ji}\frac{\partial}{\partial x}u_i \right)(x),\qquad   x\in (0,1),\; j=1,\ldots,m,
\end{equation}
where $u_j$ represents a relevant physical quantity on the $j$\textsuperscript{th} edge of the network. The operator $A$ is the gradient of the energy functional $E$ defined by
$$
E(u):= \int_0^1 \sum_{i,j=1}^m c_{ji}(x) u'_i(x)\overline{u'_j(x)} dx.
$$

As usual in the context of evolution equations on networks, we also allow for a further, weak form of interaction given by a generalized Kirchhoff-type law in the ramification nodes. These two forms of interactions between individual linear elements give rise to a well-defined system of diffusion or Schr\"odinger equations.
%
Dwelling on interesting similarities with the biological theory of neuronal coupling (cf.\ Section~5.1), we often call \emph{ephaptic} and \emph{synaptic} the influences that depend on the behaviour of the process in another edge or in another node of the network, respectively. 

Well-posedness of such diffusion and Schr\"odinger problems can be proved under quite general conditions on the coefficients $(c_{ji})$. Instead, further qualitative properties strongly depend on the coupling coefficients that are actually considered. In particular, we can show that in spite of the parabolic nature of the diffusion problem, no maximum principle holds as soon as the ephaptic coupling is nontrivial -- i.\,e., as soon as the matrix $(c_{ji})$ is non-diagonal.

\bigskip
In the second part of this paper, we discuss the issue of symmetry properties for both diffusion and Schr\"odinger equations on networks.

One says that a given physical system exhibits a \emph{symmetry} if some of its properties remain invariant under the action of a certain class of transformations. More precisely, in the Lagrangian formulation of field theory, one says that there exists a (global) symmetry of  a given dynamical system if the Lagrangian $\mathcal L(\phi)$ of the field $\phi$ is invariant under all (time- and space-independent) transformations $O$ that belong to a group $\mathcal O$, the so-called \emph{gauge group} of the system, i.\,e., if ${\mathcal L}(\phi)={\mathcal L}(O\phi)$. The prototypical example is given by the invariance under rotations of the Laplacian: this implies a symmetry for both the heat and the Schr\"odinger equations in ${\mathbb R}^n$, whose gauge group is the orthogonal group ${\mathcal O}=O_n$. Observe that since $O$ commutes with the time derivative, in many relevant cases $O$ defines a symmetry for the evolutionary problem if and only if it is a symmetry for the stationary one, i.\,e., if and only if $E(\phi)=E(O\phi)$ for all states $\phi$, where $E$ is the energy functional.

In the case of network equations, a new class of symmetries arise in a natural way: the class of proportions respected pointwise by physical quantities (e.g., temperature, densities, wave functions...) along the edges of a network during the time evolution of a physical process. To fix the ideas, consider a closed linear $Y$ subspace of ${\mathbb C}^m$ ($m$ being the number of edges in the considered network).
Then a linear closed subspace of the state space $X^2:=(L^2(0,1))^m$ can be naturally constructed as
$$\mathcal{Y} := \left\{ f \in X^2 : f(x) \in Y \mbox{ for a.\,e. } x\in (0,1) \right\}.$$
We say that $\mathcal P$ \emph{reflects a symmetry} of the network diffusion equation if the solution $u(\cdot,f)$ to the problem with initial value $f$ satisfies
\begin{equation*}
{\mathcal P} u(t,f)=u(t,{\mathcal P} f),\qquad t\geq 0,
\end{equation*}
cf.\ Definition~\ref{sym} below, where this is formulated in terms of the strongly continuous semigroup $(e^{tA})_{t \geq 0}$ generated by the operator $A$.
In Section~\ref{quantum} we will show that this is the case if and only if the orthogonal
projection $\mathcal{P}$ onto $\mathcal{Y}$ commutes with the operators of the semigroup that
governs the parabolic problem. We will also show that in the
self-adjoint case this is equivalent to the fact that
${\mathcal L}(\phi)={\mathcal L}(e^{is{\mathcal P}}\phi)$ for all $s\in\mathbb R$,
where $\mathcal L$ is the Lagrangian of the Schr\"odinger system corresponding to the
parabolic problem. In other words, we will see that \emph{${\mathcal P}$ reflects a
symmetry for the parabolic problem if and only if it generates a group of symmetries
for the Schr\"odinger system}.
In this sense $\left(e^{is\mathcal{P}}\right)_{s \in \mathbb{R}}$ can be considered as an
equivalent of a gauge group for our dynamical system. We mention that related notions of symmetries on quantum graphs have been discussed by several authors, cf.~\cite{ExnSeb89}--\cite{SevTan04}--\cite{BomKur05}.

\bigskip
Throughout this paper we will consider directed graphs. This may be disorienting at first, since we are always concerned with isotropic physical processes. In fact, all results about well-posedness as well as all those concerning positivity and asymptotics of solutions do not depend on the chosen orientation of the graph underlying the network, as it can be expected (and as it is proved in Section~3). However, we will see in Section~4 that symmetry results do in general depend on orientation: in fact, each orientation of the graph corresponds to different symmetries.

\bigskip
We will explicitely consider parabolic systems of diffusion equations in the most part of this paper. However, we will discuss in Section~\ref{quantum} how symmetry properties of both parabolic and Schr\"odinger problems can be related by means of the theory developed in Section~4, see~\cite{BolCarMug07} for more details.

\section{Well-posedness of the network equation}


The basic objects we will consider are \emph{finite directed graphs}, i.\,e., quadruples ${\mG}$ of the form $(\mV,\mE,\delta_0,\delta_1)$ where $\mV=\{\mv_1,\ldots,\mv_n\}$ and $\mE=\{\me_1,\ldots,\me_m\}$ are finite disjoint sets and $\delta_0,\delta_1:\mE\to \mV$ are mappings. They associate to an edge $\me$ two vertices $\me(0):=\delta_0(\me)$ and $\me(1):=\delta_1(\me)$, which are called \emph{initial} and \emph{terminal endpoint of $\me$}, respectively. This promptly leads to introducing two matrices $\I^+=(\iota^+_{kj})$ and $\I^-=(\iota^-_{kj})$ that fully describe the structure of the graph. They are defined by
\begin{equation}\label{incid}
\iota^+_{kj}:=\left\{
\begin{array}{rl}
1,& \hbox{if } \me _j(0)=\mv _k,\\
0, & \hbox{otherwise},
\end{array}
\right.
\quad\hbox{and}\quad
\iota^-_{kj}:=\left\{
\begin{array}{rl}
1, & \hbox{if } \me _j(1)=\mv _k,\\
0, & \hbox{otherwise.}
\end{array}
\right.
\end{equation}
Observe that if $\I:=\I^+-\I^-$ is the \emph{incidence matrix} of the directed graph $\mG$ as commonly considered in graph theory. If $|\iota_{kj}|=1$, then the edge $\me_j$ is said to be \emph{incident} to the vertex $\mv_k$. We define $$\Gamma^+(\mv _k):=\left\{j\in \{1,\ldots,m\}: \me _j(0)=\mv _k\right\}\quad\hbox{and}\quad\Gamma^-(\mv _k):=\left\{j\in \{1,\ldots,m\}: \me _j(1)=\mv _k\right\},$$
and by $\Gamma(\mv_k):=\Gamma^+(\mv_k) \cup \Gamma^-(\mv_k)$ we denote the set of indices of all edges that are incident to $\mv_k$.

If there exists an edge $\me\in\mE$ such that either $\me(0)=\mv_k$ and $\me(1)=\mv_\ell$, or $\me(0)=\mv_\ell$ and $\me(1)=\mv_k$, then the vertices $\mv_k,\mv_\ell$ are said to be \emph{adjacent}. Similarly, we say that edges $\me_i,\me_j$ are \emph{adjacent} if there exists a vertex which they are both incident to, i.\,e., if there exists $\mv\in \mV$ such that $\me_i(0)=\mv$ or $\me_i(1)=\mv$, and such that $\me_j(0)=\mv$ or $\me_j(1)=\mv$.


\bigskip
Additionally, we assign to the graph a metric structure that allows us to treat it as a one-dimensional manifold and, eventually, to consider partial differential equations describing evolution processes taking place on it. Throughout this paper we will always call \emph{network} any directed graph endowed with such a metric structure. A similar if not identical approach, based on von Below's theory of $C^2$-networks, has been presented in~\cite{Bel85}.

More precisely, each edge of the graph will be thought of as an interval. For the sake of consistency with the notation introduced in~\eqref{incid}, such intervals are parametrized in such a way that they have length $1$. Whenever we consider a square integrable function $f$ acting on the graph $\mG$, we may equivalently think of $f$ as a complex-valued function $\mG\to{\mathbb C}$ defined almost everywhere (with respect to the $1$-dimensional Lebesgue measure) on the edges of the graph, or equivalently as a vector-valued function $(0,1)\to{\mathbb C}^m$. In this case we will denote $f$ by $(f_1,\ldots,f_m)^\top$, where each $f_j\in L^2(0,1)$ is a function on $\me_j$, $i=1,\ldots,m$. Whenever point evaluations of $f$ are well-defined, we define with an abuse of notation $f_j(\mv_k):=f_j(0)$ if $\iota^+_{kj}=1$, and $f_j(\mv_k):=f_j(1)$ if $\iota^-_{kj}=1$.

\bigskip
As already emphasized in Section~1, in contrast with the setting which is usual in the literature on network evolution equations, we discuss a general model and allow for (possibly non-mutual) interactions of non-adjacent pairs of edges, too. The influence of the process taking place along the edge $\me_i$ onto that taking place along $\me_j$ will be descrived by the \emph{ephaptic} coupling coefficient $c_{ji}$. Such a coefficient is seen as a function on the edge $\me_i$: with the same convention as above we thus denote $c_{ji}(\mv_\ell):=c_{ji}(0)$ or $c_{ji}(\mv_\ell):=c_{ji}(1)$ if $\iota ^+_{\ell i}=1$ or $\iota^-_{\ell i}=1$, respectively.

While the dynamics of our system is described by the coupled diffusion equations in \eqref{parab},  we still have to equip it with suitable conditions in the nodes. To this aim, we introduce two tensors defined by 
$${\mathfrak I}^+ := \I^+ \otimes \I^+\qquad\hbox{and}\qquad {\mathfrak I}^-:= \I^- \otimes \I^-.$$ 
We call $\mathfrak I:=\mathfrak I^+ -  {\mathfrak I}^-$ the \emph{ephaptic incidence tensor} of $\mG$. Here $\otimes$ stands for the usual Kronecker product of two $m\times n$ matrices, defined by $(A\otimes B)^{kj}_{\ell i}:=a_{kj}\cdot b_{\ell i}$.
We denote by ${\iota}^{kj}_{\ell i},\hat{\iota}^{kj}_{\ell i},\check{\iota}^{kj}_{\ell i}$ the entries of $\mathfrak I,\mathfrak I_+,\mathfrak I_-$, respectively.
In other words, $\iota^{kj}_{\ell i}$ represents the influence of the vertex $\mv_\ell$ as an endpoint of $\me_i$ on the vertex $\mv_k$ as an endpoint of $\me_j$. By construction, such influences are symmetric, i.\,e., ${\iota}^{kj}_{\ell i}={\iota}_{kj}^{\ell i}$ for all $i,j=1,\ldots,m$ and all $k,\ell=1,\ldots,n$. 

\bigskip
Solutions of our network diffusion problem have to be continuous in the vertices, i.\,e.,
\begin{equation}\label{contin}
u_i(\mv_k)= u_j(\mv_k)\qquad \hbox{for all }i,j\in \Gamma(\mv_k),\; k=1,\ldots,n.
\end{equation}
Because of the continuity condition expressed in the equation \eqref{contin}, we can and  will denote by $d^u_k$ the joint value of the components of the vector-valued function $u$ at the node $\mv_k$.

Furthermore, we allow (possibly non-adjacent) vertices of the graph to influence each other. A natural interaction condition can be formulated as
\begin{equation*}\label{epsyn}
\sum\limits_{i,j=1}^m \sum_{k=1}^n \omega_{\ell i}^{kj} u'_i(t,\mv_\ell)=0,\qquad k=1,\ldots,m.
\end{equation*}
Here, the {weighted} incidence tensor $\mathfrak W:=(\omega_{\ell i}^{kj})$, for $i,j=1,\ldots,m$ and $\ell,k=1, \ldots, n$, is defined by
\begin{equation*}\label{weitensor}
\omega_{\ell i}^{kj}:=c_{ji}(\mv_\ell)\iota_{\ell i}^{kj}.
\end{equation*}

In fact, in a fashion similar to that considered in~\cite{Mug07a} we allow for even more general, non-local Kirchhoff-type conditions. Such generalized conditions are given by
\begin{equation}\label{kkk}
\sum\limits_{i,j=1}^m \sum\limits_{\ell=1}^n \omega_{\ell i}^{kj} u'_i(t,\mv_\ell)=\sum\limits_{\ell=1}^n m_{k\ell }d^u_\ell, \qquad k=1,\ldots,m.
\end{equation}

Summing up, we investigate the strongly coupled system of initial-boundary value diffusion problems
\begin{equation}\label{problem}
\left\{
\begin{array}{rcll}
\dot{u}_j(t,x)&=&\sum\limits_{i=1}^m (c_{ji}u'_i)'(t,x),&x\in(0,1),t>0,\;j=1,\ldots,m,\\
u_i(t,\mv_k)&=& u_j(t,\mv_k)=: d^u_k(t), &t\geq 0,\; i,j\in \Gamma(\mv_k),\; k=1,\ldots,n,\\
\sum\limits_{\ell=1}^n m_{k\ell}d^u_\ell &=&\sum\limits_{i,j=1}^m \sum\limits_{\ell=1}^n \omega_{\ell i}^{kj}u'_j(t,\mv_\ell),& t\geq 0,\; k=1,\ldots,n,
\\
u_j(0,x)&=&u_{j0}(x), &x\in(0,1),\;j=1,\ldots,m.
\end{array}
\right.
\end{equation}


\begin{rem}\label{kirchgen}
By definition of ${\mathfrak W}$, whenever $C(x)\equiv \Id$ and $M=0$ (i.\,e., if only local, synaptic interaction occurs),~\eqref{kkk} reduces to the usual Kirchhoff condition prescribing that in each node $\mv_\ell$ incoming and outgoing heat fluxes agree.
\end{rem}

We introduce $X^2:=(L^2(0,1))^m$, which is a Hilbert space with respect to the canonical inner product
$$( f\mid g)_H=\sum_{j=1}^m \int_0^1 f_j(x) \overline{g_j(x)}dx,\qquad f,g\in V.$$
We also consider its dense subspace
$$V:=\{f\in(H^1(0,1))^m: \exists d^f\in \C^n \mbox{ s.\,t. } (\I^+)^\top d^f=f(1),(\I^-)^\top d^f=f(0)\},$$
the space of all $H^1$-functions that are continuous in the nodes of the graph. The subspace $V$ is a Hilbert space with respect to the canonical inner product
$$(f\mid g)_V:=\sum_{i=1}^m \int_0^1 \left( f'_i(x)\overline{g'_i(x)} +f_i(x)\overline{g_i(x)}\right)dx,\qquad f,g\in V.$$
Observe that $V$ is densely and compactly imbedded into $X^2$, since $(C_c^\infty(0,1))^m \subset V\subset (L^2(0,1))^m$. 


\bigskip
For the sake of later reference, we recall that a complex (possibly nonsymmetric) matrix $M=(m_{ij})$ is called  \emph{accretive} (resp., \emph{dissipative}) if there exists $\mu\geq 0$ such that ${\rm Re}(M\xi| \xi)\geq \mu |\xi|^2$ (resp., ${\rm Re}(M\xi| \xi)\leq -\mu|\xi|^2$)
for all $\xi\in{\mathbb C}^m$. We call $M$ \emph{positive definite} (resp., \emph{negative definite}) if it is accretive (resp., dissipative) and moreover $\mu$ can be chosen $>0$.

Throughout the remainder of this paper we will always assume the following.

\begin{assum}\label{coercassum1}
The coefficients $c_{ij}$ are functions of class $C^1[0,1]$. The matrix $C(x)=(c_{ij}(x))$ is positive definite, uniformly on the interval $[0,1]$, i.\,e., there exists $\mu>0$ such that
\begin{equation*}
\Real (C(x)v\mid v)_{{\mathbb C}^m}:=\Real \sum_{i,j=1}^m c_{ij}(x)v_j \overline{v_i} \geq \mu |v|^2_{{\mathbb C}^m} \qquad\hbox{for all } x\in [0,1],\; v\in{\mathbb C}^m.
\end{equation*}
\end{assum}

Observe that Assumption~\ref{coercassum1} is weaker than~\cite[Assum.~2.3]{CarMug07}.

%
%

Let us now introduce the Kirchhoff operators $\Phi^+,\Phi^-:(H^2(0,1))^m\to {\mathbb C}^n$ defined by
\begin{equation*}
\Phi^+u:=\begin{pmatrix}
\sum\limits_{i,j=1}^m\sum\limits_{\ell=1}^n \hat{\omega}^{j1}_{\ell i}u'_i(\mv_1) \\
\vdots\\
\sum\limits_{i,j=1}^m\sum\limits_{\ell=1}^n\hat{\omega}^{jn}_{\ell i}u'_i(\mv_n) \\
        \end{pmatrix},\qquad 
\Phi^-u:=\begin{pmatrix}
\sum\limits_{i,j=1}^m\sum\limits_{\ell=1}^n \check{\omega}^{j1}_{\ell i}u'_i(\mv_1) \\
\vdots\\
\sum\limits_{i,j=1}^m\sum\limits_{\ell=1}^n \check{\omega}^{j1}_{\ell i} u'_i(\mv_n) \\
        \end{pmatrix},
\end{equation*}
and a differential operator on $X^2$ by
\begin{equation}\label{operator}
A:=\left(
\begin{array}{cccc}
\frac{d}{dx}(c_{11}\frac{d}{dx}) &  \dots & \frac{d}{dx}(c_{1m}\frac{d}{dx}) \\
\vdots&\ddots&\vdots\\
\frac{d}{dx}(c_{m1}\frac{d}{dx}) & \dots & \frac{d}{dx}(c_{mm}\frac{d}{dx}) \\
\end{array}
\right)
\end{equation}
with domain
\begin{equation}\label{domain}
D(A):=\left\{f\in (H^2(0,1))^m \cap V: \Phi^+f-\Phi^-f= Md^f \right\},
\end{equation}
for the matrix $M=(m_{kh})$ introduced in~\eqref{kkk}. Since $D(A)\subset V$, functions in $D(A)$ are continuous in the nodes.

With the aim of pursuing a variational approach to our problem, we introduce a densely defined sesquilinear form $a$ defined by
\begin{equation}\label{form}
a (f,g):=(Cf' \mid g')_{X^2}-(Md^f \mid d^g)_{\mathbb C^n}=\sum\limits_{i,j=1}^m \int\limits_0^1c_{ij}(x)f'_j(x)\overline{g'_i(x)}dx - \sum\limits_{k,\ell=1}^n m_{k\ell}d^f_\ell\overline{d^g_k}
\end{equation}
for $f,g\in V$, which will be later shown to be related to the operator $A$.

%

\begin{theo}\label{wp}
The operator associated with the form $a$ generates a compact, analytic semigroup on $X^2$. This semigroup is contractive (hence asymptotically almost periodic, too) if  $M$ is dissipative. If $M$ is dissipative, then the semigroup is strongly stable if and only if $M^*{\mathbb 1}\not=0$. The semigroup is uniformly exponentially stable if $M$ is negative definite. The semigroup is self-adjoint if and only if the matrices $C(x)$, $x\in[0,1]$, and $M$ are self-adjoint.
\end{theo}

Observe that the last result also characterizes well-posedness of the quantum graph associated with~\eqref{problem}.

We stress that if the semigroup is contractive (resp., uniformly exponentially stable), then $M$ is not necessarily dissipative (resp., negative definite), as one sees already in the case of a network consisting of a single interval, if one considers the function $f$ defined by $f(x) = x$ and $M=\Id$.

\begin{proof}
We show that the sesquilinear form $ a$ is continuous and $X^2$-elliptic, i.\,e.,
\begin{itemize}
\item $| a (f,g)|\leq K_1 \left\|f\right\|_V \left\|g\right\|_V$ for some constant $K_1>0$ and all $f,g\in V$, and
\item there exist $\alpha>0$ and $\omega\in\mathbb R$ such that $\Real  a (f,f)\geq \alpha \Vert f\Vert^2_{V}-\omega \| f\|^2_{X^2}$ for all $f\in V,$
\end{itemize}
respectively. In fact, the continuity of $ a $ is a direct consequence of the Cauchy--Schwarz inequality in $X^2$ and of the continuous imbedding of $V$ into $(C[0,1])^m$, and the constant $K_1$ is the maximum over $x\in[0,1]$ of the matrix norm $\| C(x)\|$.

In order to prove $X^2$-ellipticity of $a$, it suffices to observe that $(Cf'\mid g')_{X^2}$ clearly defines an $X^2$-elliptic form if (and only if) $C(x)$ is a positive definite matrix for a.\,e. $x\in[0,1]$, which is Assumption~\ref{coercassum1}. Since there exists $K_2>0$ such that
$$\max_{x\in[0,1]} |f(x)|\leq K_2 \| f\|^{\frac{1}{2}}_{L^2} \|f\|^{\frac{1}{2}}_{H^1},\qquad f\in H^1(0,1),$$
cf.~\cite[Cor.~4.11]{Bur98}, it follows that the space of continuous functions over the graph is an interpolation space between $(H^1(0,1))^m$ and $(L^2(0,1))^m$. It then suffices to apply~\cite[Lemma~2.1]{Mug07b} in order to treat the lower order perturbation given by $(Md^f\mid d^g)_{{\mathbb C}^n}$. Accordingly, by~\cite[Prop.~1.51 and Thm.~1.52]{Ouh04} the operator associated with $a$ generates an analytic semigroup of angle $\frac{\pi}{2}-\arctan K_1$.

Observe that by the Rellich--Kondrachov theorem the embedding of $V$ into $X^2$ is compact, thus the semigroup is compact. A direct computation shows  that $a$ is accretive (i.\,e., $\Real  a (f,f)\geq 0$ for all $f\in V$) if $M$ is dissipative; and that $ a $ is coercive (i.\,e., it is $X^2$-elliptic with $\omega=0$) if $M$ is negative definite.
In the first case, the semigroup associated with $a$ is contractive, and by~\cite[Thm.~5.5.6]{AreBatHie01} also asymptotically almost periodic. In the latter case, the semigroup is uniformly exponentially stable since the shifted form $a-\alpha(\cdot |\cdot)_V$ is accretive. Finally, let $M$ be dissipative. Then, by~\cite[Exa.~V.2.23]{EngNag00} the semigroup associated with $a$ is strongly stable if and only if $0$ is not an eigenvalue of the operator associated with the adjoint form $a^*$. First of all, observe that if $A^*f =0$, then necessarily
$$\mu \| f'\|^2_{X^2}\leq (Cf'|f')_{X^2}=(Md^f | d^f)_{{\mathbb C}^n}\leq 0,$$ 
thus $f$ is a constant, i.\,e., $f=c\mathbb 1$. Observe now that $0$ is an eigenvalue of $A^*$ (and thus necessarily with eigenfunction $\mathbb 1$) if and only if $$0=(A^*{\mathbb 1}|g)=-a^*({\mathbb 1},g)=(M^* {\mathbb 1} | d^g)_{{\mathbb C}^n}\qquad \hbox{for all }g\in V,$$
and since the nodal values $d^g$ of $g$ are arbitrary vectors of ${\mathbb C}^n$, this is equivalent to saying that $M^*{\mathbb 1}=0$. Finally,  $a$ is self-adjoint if and only if so are the coefficient matrices.
\end{proof}

\begin{rem}\label{ahiahiahi}
It is known that the operator associated with the form $a$ cannot generate an analytic, quasicontractive semigroup unless $a$ is $X^2$-elliptic, (see~\cite[\S~5.3.4]{Are04}), and hence unless Assumption~\ref{coercassum1} holds.
\end{rem}

In order to show the well-posedness of our motivating problem, we need to make sure the operator associated with $a$ is actually $A$ as introduced in~\eqref{operator}--\eqref{domain}. Having proved this, Theorem~\ref{wp} becomes a generation result for $A$, and in the remainder of this paper we will denote by $(e^{tA})_{t\geq 0}$ the semigroup introduced above.

\begin{prop}
The operator associated with $ a $ is $(A,D(A))$ as defined in \eqref{operator}--\eqref{domain}.
\end{prop}

\begin{proof}
Denote by $(B,D(B))$ the operator associated with the form $ a $, which by definition is given by
\begin{equation*}
\begin{array}{rcl}
D(B)&:=&\left\{f \in V: \exists g \in X^2 \mbox{ s.\,t. }  a (f,h)=(g \mid h)_H\; \forall h\in V\right\},\\
Bf&:=&-g.
\end{array}
\end{equation*}
We first show that $A\subset B$. Fix $f\in D(A)$. Then for all $h \in V$
\begin{eqnarray*}
 a (f,h) &=&\sum\limits_{i,j=1}^m \int\limits_0^1c_{ji}(x)f'_i(x)\overline{h'_j(x)}dx - \sum\limits_{k,\ell=1}^n m_{k\ell}d^f_k\overline{d^h_\ell}\\
&=&\sum\limits_{i,j=1}^m [c_{ji}f'_i\overline{v_j}]^1_0 - \sum\limits_{i,j=1}^m \int\limits_0^1(c_{ji}f'_i)'(x)\overline{h_j(x)}dx-\sum\limits_{k,\ell=1}^n m_{k\ell}d^f_\ell\overline{d^h_k}.
\end{eqnarray*}
Using now the definition of the incidence tensor $\mathfrak I=\hat{\mathfrak I}-\check{\mathfrak I}$ we can write
\begin{eqnarray*}
\sum\limits_{i,j=1}^m [c_{ji}f'_i\overline{h_j}]^1_0 &=&\sum\limits_{i,j=1}^m\sum\limits_{\ell, k=1}^n c_{ji}(\mv_\ell) (\hat{\iota}^{kj}_{\ell i}-\check{\iota}^{kj}_{\ell i}) f'_i(\mv_\ell)\overline{h_j(\mv_k)}\\
&=&\sum\limits_{i,j=1}^m\sum\limits_{k,\ell=1}^n {\omega}^{kj}_{\ell i} f'_i(\mv_\ell)\overline{h_j(\mv_k)}\\
&=&\sum_{k=1}^n \overline{d^h_k}\sum_{i,j=1}^m\sum_{\ell=1}^n \omega^{kj}_{\ell i} f'_i(\mv_\ell)\\
&=&\sum_{\ell,k=1}^n m_{k \ell} d^f_\ell \overline{d^h_k}.
\end{eqnarray*}
As a consequence
$$
a(f,h)=- \sum\limits_{i,j=1}^m \int\limits_0^1(c_{ji}f'_i)'(x)\overline{h_j(x)}dx=: ((Cf')' \mid h).
$$
Thus, for all $h\in V$ there exists $g=Af\in X^2$ such that
\begin{equation*}
 a (f,h)=- \sum\limits_{i,j=1}^m \int\limits_0^1(c_{ji}f'_i)'(x)\overline{h_j(x)}dx=-(g \mid h)_H.
\end{equation*}
This completes the proof of the first inclusion. Conversely, let $f \in D(B)$. By definition there exists $g \in V$ such that $ a (f,h)=-(g\mid h)_H$ for all $h \in V$, and accordingly
\begin{equation*}
\sum\limits_{i,j=1}^m \int\limits_0^1c_{ji}(x)f'_i(x)\overline{h'_j(x)}dx - \sum\limits_{\ell, k=1}^n m_{\ell k}d^f_k\overline{d^h_\ell}=- \sum\limits_{i=1}^m \int\limits_0^1g_i(x)\overline{h_i(x)}dx.
 \end{equation*}
Integrating by part the left hand side, we obtain that
$$ -  \sum\limits_{i,j=1}^m \int\limits_0^1(c_{ji}f'_i)'(x)\overline{h_j(x)}dx 
 + \sum\limits_{k=1}^n \overline{d^h_k}\sum\limits_{i,j=1}^m\sum\limits_{\ell=1}^n \omega^{kj}_{\ell i} f'_i(\mv_\ell) 
- \sum\limits_{ k,\ell=1}^n m_{k\ell}d^f_\ell\overline{d^h_k} =
-   \sum\limits_{i=1}^m \int\limits_0^1g_i(x)\overline{h_i(x)}dx,$$
which holds for all $h \in V$. In particular, considering $h\in (H^1_0(0,1))^m$ 
vanishing on all but one edge of the network, we conclude that 
\begin{equation*}\label{coincideriv}
g_{i}(x)=\sum\limits_{j=1}^m(c_{j{i}}f_j')'(x) \qquad \hbox{for all } x\in (0,1) \hbox{ and all }i=1,\ldots,m.
\end{equation*}
Similarly, considering $h$ with arbitrary nodal values and arbitrary small $X^2$-norm, we obtain
\begin{equation*}
\sum\limits_{i,j=1}^m\sum\limits_{\ell=1}^n \omega^{kj}_{\ell i} f'_i(\mv_\ell) 
- \sum\limits_{\ell=1}^n m_{k\ell}d^f_\ell = 0\qquad\hbox{for all } k=1,\ldots,n.
\end{equation*}
This shows that $f \in D(A)$ and completes the proof.
\end{proof}


Having proved analytical well-posedness in an $L^2$-space, one could try to extend this result to further $L^p$-spaces, $p\not= 2$. To this end, a common strategy is to show that the semigroup leaves invariant the unit ball of $L^\infty$, so that each operator $e^{tA}$, $t\geq 0$ is contractive on all $L^p$ spaces, $p\in[2,\infty]$, by virtue of Riesz--Thorin interpolation theorem. This has already been accomplished in the case of pure synaptic coupling, cf. \cite{KraMugSik07}--\cite{Mug07a}. However, we show in the following that this approach cannot work in the case of nontrivial ephaptic coupling.

\begin{theo}\label{contractive}
The following assertions hold.
\begin{enumerate}
\item The semigroup $(e^{tA})_{t\geq 0}$ is real, i.\,e., it leaves invariant the subspace of real-valued function of $X^2$, if and only if 
\begin{itemize}
\item $C(x)\in M_m({\mathbb R})$  for all $x\in [0,1]$ and 
\item $M\in M_n({\mathbb R})$.
\end{itemize}
\item The semigroup $(e^{tA})_{t\geq 0}$ is positive, i.\,e., it leaves invariant the positive cone of $X^2$, if and only if 
\begin{itemize}
\item $C(x)$ is a real valued, diagonal matrix  for all $x\in [0,1]$ and 
\item the matrix $M$ has real entries that are positive off-diagonal.
\end{itemize}
In this case, the semigroup is also irreducible if the graph is connected.
\item The semigroup $(e^{tA})_{t\geq 0}$ is $X^\infty$-contractive (resp., $X^1$-contractive), i.\,e., it leaves invariant the unit ball of $X^\infty$ (resp., of $X^1$), if and only if
\begin{itemize}
\item $C(x)$ is a real valued, diagonal matrix  for all $x\in [0,1]$ and 
\item the matrix $M$ satisfies 
$\Real  m_{kk} + \sum_{h\not=k} \vert m_{kh}\vert \leq 0$ (resp., $\Real  m_{kk} + \sum_{h\not=k} \vert m_{hk}\vert \leq 0$) for all $k=1,\ldots,n.$
\end{itemize}
\end{enumerate}
\end{theo}

\begin{proof}
As shown in the proof of Theorem~\ref{wp} the form  $ a $ is densely defined, continuous, and $X^2$-elliptic. Thus, by~\cite[Prop. 2.5]{Ouh04}, and by a simple rescaling argument, the semigroup $(e^{t a })_{t\geq 0}$ is real if and only if $\Real f\in V$ and $ { a }(\Real {f},\Range f)\in\mathbb{R}$ for all $f\in V$. Thus, an easy computation shows that reality of the coefficients $C,M$ is sufficient.

Conversely, assume $(e^{tA})_{t\geq 0}$ to be real. Let $f_0\in H^1_0(0,1)$ real valued and such that its support of $f_j$ agrees with $[a,b]\subset (0,1)$. Define $f$ as a function such that $f_i=if_0$, $f_j=f_0$, and all further coordinates vanish. By the above characterization of real semigroups one has $a(\Real {f},\Range f)=\int_a^b c_{ij}(x) |f'_0(x)|^2 dx\in\mathbb{R}$. Since this construction can be repeated for arbitrary $a,b$ and $i,j$, we deduce that $c_{ij}(x)$ is a real number  for all $x\in(0,1)$, and by continuity also for all $x\in [0,1]$.

Let now $f\in V$ such that $d^f_{\ell}=1$ and $d^f_k=i$. If $f$ vanishes in all further nodes, $a(\Real {f},\Range f)=(C (\Real f)'\mid (\Range  f)')_{X^2}-m_{k\ell}$. As shown above, $C(x)$ is a real matrix for all $x\in [0,1]$ and therefore $(C (\Real f)'\mid (\Range  f)')_{X^2}\in\mathbb R$. Thus, $m_{k\ell}\in\mathbb R$ for all $k,\ell=1,\ldots,n$.

In a similar fashion and taking into account~\cite[Thm.~3.5]{Mug07a} and~\cite[Prop.~3.6]{CarMug07}, one can prove the claimed characterizations of positivity, $X^\infty$-contractivity and, by duality,  $X^1$-contractivity of $(e^{tA})_{t\geq 0}$.
\end{proof}

Additional properties of boundary regularity of solutions of~\eqref{problem} can be deduced by the fact that the analytic semigroup operators $e^{tA}$ map $X^2$ into $\bigcup_{k=1}^\infty D(A^k)$ for all $t>0$.

\begin{prop}
If $u$ is the solution to~\eqref{problem}, the following assertions hold.
\begin{enumerate}
\item $\sum_{j=1}^m (c_{ij} u_j')'$ is continuous in the nodes and satisfies a Kirchhoff law, i.\,e.,
\begin{eqnarray*}
\sum_{\iota=1}^m (c_{i\iota} u_\iota')'(t,\mv_\ell)&=& \sum_{\iota=1}^m (c_{j\iota } u_\iota')'(t,\mv_\ell)=:d^{(cu')'}_\ell(t),\qquad t> 0,\; i,j \in \Gamma(\mv_\ell),\; \ell=1,\ldots,n,\\
\sum\limits_{\ell=1}^n m_{k\ell}d^{(cu')'}_\ell(t) &=&\sum\limits_{\iota,i,j=1}^m \sum\limits_{\ell=1}^n \omega_{\ell \iota}^{kj}(c_{j \iota }u_\iota ')''(t,\mv_\ell),\qquad t> 0,\; k=1,\ldots,n.
\end{eqnarray*}
\item If furthermore the coefficients matrix $C$ is diagonal, then $u$ is of class $C^\infty$ and its derivatives of even and odd order satisfy for all $N\in\mathbb N$
\begin{eqnarray*}
u^{(2N)}_i(t,\mv_\ell)&=& u^{(2N)}_j(t,\mv_\ell)=: d^{u^{(2N)}}_\ell(t),\qquad t> 0,\; i,j\in \Gamma(\mv_\ell),\; \ell=1,\ldots,n,\\
\sum\limits_{\ell=1}^n m_{k\ell}d^{u^{(2N)}}_\ell(t) &=&\sum\limits_{i,j=1}^m \sum\limits_{\ell=1}^n \omega_{\ell \iota}^{kj}u^{(2N+1)}_j(t,\mv_\ell),\qquad t\geq 0,\; k=1,\ldots,n.
\end{eqnarray*}
\end{enumerate}
\end{prop}

\section{Symmetry Properties}

In this section we will characterize invariance of different classes of closed linear subspaces of the space $X^2$ under the action of $(e^{tA})_{t\geq 0}$. The invariance of a closed subspace under the action of a semigroup can be characterized as a direct consequence of a result due to E.-M.\ Ouhabaz, see~\cite[Thm.~2.2]{Ouh04}. For the sake of self-containedness we present it in the form we will use in the following. Observe that in the view of~\cite[Cor.~5.2]{CarMug07}, the invariance results for susbpaces deduced by means of Theorem~\ref{oulin} can be directly extended to a large class of nonlinear, strip-like subsets of $X^2$.

\begin{theo}\label{oulin}
Let $a:V\times V\to\mathbb C$ be a continuous, elliptic sesquilinear form on a Hilbert space $H$, and consider an orthogonal projection $P$ on $H$. Then $\Range \mathcal P$ is invariant under the action of the semigroup $(e^{tA})_{t\geq 0}$ associated with $a$ if and only if
\begin{enumerate}
\item $PV\subset V$ and
\item $a(f,g)=0$ for all $f \in \Range \mathcal P \cap V, g\in \Ker {\mathcal P} \cap V$.
\end{enumerate}
\end{theo}

%
%
%
%
%
A relevant class of subspaces of $X^2$ can be constructed as follows: Let $Y$ be a subspace of $\mathbb{C}^m$ and consider
\begin{equation}\label{projy}
\mathcal{Y} := \left\{ f \in X^2 : f(x) \in Y \mbox{ for a.\,e. } x\in (0,1) \right\}.
\end{equation}
We look for criteria for invariance of the subspace $\mathcal Y$ of $X^2$ under the action of the semigroup $(e^{tA})_{t\geq 0}$. Denoting by $K$ the orthogonal projection of ${\mathbb C}^m$ onto $Y$, the orthogonal projection $\mathcal{P}_K$ of $X^2$ onto $\mathcal{Y}$
satisfies
\begin{equation}\label{projdef}
\left(\mathcal{P}_Kf\right)(x) = K\left(f(x)\right)\qquad \hbox{for a.\,e. } x\in(0,1).
\end{equation}

\bigskip
The aim of this section is to discuss problems that are similar to that presented in the following, which also shows
an intuitive relation between invariance and symmetry properties.

\begin{exa}\label{motivexa}
Consider a graph $\mG$ consisting of two edges, both outgoing from a common vertex $\mv$, i.\,e., an outbound star. Let $C=\Id$ and $M=0$. Then the form $a$ is associated with the Laplacian with a Kirchhoff condition in $\mv_1$ and Neumann conditions in the boundary nodes. \emph{Do initial data that are symmetric with respect to $\mv$ give rise to solutions to the diffusion problem that are also symmetric with respect to $\mv$?} We can reformulate this question and ask whether the closed linear subspace $\mathcal Y:=\{f \in X^2: f_1=f_2\}$ is invariant under the action of the semigroup $(e^{tA})_{t\geq 0}$. In fact, $\mathcal Y=\Range \mathcal {\mathcal P}_K$, where $K$ is the $2\times 2$ matrix whose entries equal $\frac{1}{2}$.
\end{exa}

Let us reformulate the criterion in Theorem~\ref{oulin} in our special case. After rewriting the form $a$ as $
a(f,g)=(Cf'\mid g')_{X^2} - (Md^f \mid d^g)_{{\mathbb C}^n}$,  observe that the denseness of $V_{x}:=\{f \in V: d^f = x\}$ in $X^2$ for each $x\in \mathbb C^n$ implies that the condition (2) of Theorem \ref{oulin} holds if and only if 
\begin{equation}\label{Cortho}
(Cf'\mid g')_{X^2}=0\qquad \mbox{ for all }f \in \Range {\mathcal P} \cap V, g\in \Ker {\mathcal P} \cap V
\end{equation}
and 
\begin{equation}\label{Mortho1}
(Md^f \mid d^g)_{\mathbb C^n}=0 \qquad \hbox{for all }f \in \Range \mathcal {\mathcal P} \cap V, g \in \Ker \mathcal {\mathcal P} \cap V.
\end{equation}
We will refer to condition (1) of Theorem \ref{oulin} as to the \emph{admissibility of the projection $\mathcal {\mathcal P}_K$} (or sometimes of $K$), and to the condition \eqref{Cortho} and \eqref{Mortho1} as the \emph{orthogonality condition with respect to $\mathcal {\mathcal P}_K$} of the coefficient matrix $C$ and of the matrix $M$, respectively. Characterizing admissibility and orthogonality is aim of the following subsections.
\subsection{Admissibility}
In particular, $ \Ker \mathcal {\mathcal P}_K$ and $\Range \mathcal {\mathcal P}_K$ are isomorphic to $(L^2(0,1))^k$ and $(L^2(0,1))^r$, respectively.

We will now investigate the admissibility of projections of the type $\mathcal {\mathcal P}_K$ in terms of the matrix $K$ and of (the incidence matrix $\I$ of) the graph $\mG$. 
Let us fix some notation. For $\mathcal A \subset \{1,\ldots,m\}$ we define the vector
\begin{equation}\label{indicator}
\1_\mathcal A:=(a_i)_{i=1,\ldots,m}, \mbox{ where } a_i:= \begin{cases} 1 &i\in\mathcal A,\\0 & i\not\in \mathcal A\end{cases}
\end{equation}
and write $\mathbb{1} := \mathbb{1}_\mathcal A$ in the special case of $\mathcal{A} = \{1,\ldots,m\}$.

\begin{lemma}\label{1eigen}
Let the graph $\mG$ be connected and the projection $\mathcal {\mathcal P}_K$ be admissible. Then $\1$ is an eigenvector of $K$.
\end{lemma}
\begin{proof}
By hypothesis $\mathcal {\mathcal P}_K V\subset V.$ Consider the function $\1: x \mapsto (1,\ldots,1)^\top$ and observe that $\mathcal {\mathcal P}_K \1(x)=K\1$ and $\mathbb 1 \in V$. This shows that on each edge $\mathcal {\mathcal P}_K \1$ is a constant function, and since $\mathcal {\mathcal P}_K \1 \in V$ all these constants coincide, hence $\mathcal {\mathcal P}_K \1 = a\1$ for an $a \in \mathbb C$.
\end{proof}
\begin{rem}
Observe that $K\1 \in \{0,\1\}$, since the only eigenvalues of an orthogonal projection are $0$ and $1$, and that $\1 \in \Ker (\Id-K)$ if $\1 \in \Range K$. Moreover, $K$ is admissible if and only if $\Id-K$ is admissible. Therefore we may assume $\1 \in \Range K$ without loss of generality.
\end{rem}

Lemma \ref{1eigen} can be used to investigate the invariance of subgraphs.

\begin{exa}
If the graph $\mG$ is connected, then there exists no proper subgraph $\mG'$ of $\mG$ such that the linear subspace  $\mathcal{Y}:=\{f \in X^2: f_{|G'}=0\}$ of the functions vanishing on $\mG'$ is invariant under the action of $(e^{tA})_{t\geq0}$.

Without loss of gererality we may assume that the subgraph $\mG'$ corresponds to the edges $\me_{m'+1},\ldots, \me_m$. The projection onto $Y$ is given by $\mathcal {\mathcal P}_K$, where $$K=\begin{pmatrix}\Id_{m'}&0\\0&0\end{pmatrix}.$$ 
Of course, $\1$ is not an eigenvector of $K$. This result is independent of the matrices $C$ and $M$.
\end{exa}

To characterize admissibility of projections having $\1$ as an eigenvector we introduce some additional notation. We define the $2m\times n$ matrix $\tilde{\mathcal I}$ and the $2m\times 2m$ matrix $\tilde{K}$ as
\begin{equation}\label{deftilde}
\It:=(\I^+,\I^-)^\top=\begin{pmatrix} (\I^+)^\top \\(\I^-)^\top\end{pmatrix} \qquad \text{and} \qquad \K:=\begin{pmatrix} K&0\\0&K\end{pmatrix}.
\end{equation}
Observe that $\K$ is an orthogonal projection of $\mathbb C^{2m}$.

\begin{lemma}\label{decomp}
Let the matrix $K$ be an orthogonal projection of $\mathbb C^d$ and the let the set $Y$ be a linear subspace of $\mathbb C^d$. Then the following assertions are equivalent.
\begin{enumerate}[(a)]
\item $KY \subset Y$;
\item $Y=\Ker K \cap Y \oplus \Range  K \cap Y$.
\end{enumerate}
\end{lemma}
\begin{proof}
``(b) $\Rightarrow$ (a)''. Let $u \in Y$, i.\,e., $u=u_1 +u_2$, where  $u_1 \in \Ker K \cap Y$ and $u_2 \in \Range  K \cap Y$. Then $Ku=Ku_1 + Ku_2=u_2 \in Y$, which proves the claim.\\
``(a) $\Rightarrow$ (b)''. Let $B^1=\{b^1_i: i=1,\ldots,r_0\}$ be a basis of $\Ker  K \cap Y$ and $B^2=\{b^2_i: i=1,\ldots,q_0\}$ be a basis of $\Range  K \cap Y$. Extend $B^1$ and $B^2$ to a basis of $\Ker K$ and $\Range K$, respectively, denoted by
$$B^{1\star}=B^1 \cup \{b^1_i: i=r_0+1,\ldots,r\},$$
and
$$B^{2\star}=B^2 \cup \{b^2_j: j=q_0+1,\ldots,q\}.$$
Observe that $\mathbb C^d=\Ker K \oplus \Range K$ since $K$ is an orthogonal projection. Let $u \in Y$. Then
$$
u=\sum_{i=1}^{r} \alpha_i b^1_i + \sum_{i=1}^{q} \beta_i b^2_i
$$
with uniquely determined coefficients $\alpha_i,\beta_j, i=1,\ldots,r, j=1,\ldots,q.$ Now
$$
Ku- \sum_{i=1}^{q_0} \beta_i b^2_i = \sum_{i=q_0+1}^{q} \beta_i b_i^2 \in \Range K \cap Y,
$$
since $Ku \in Y$ by assumption, and hence $\beta_j=0, j=q_0+1,\ldots,q$ by definition of $B^{2\star}$.
Analoguously it can be shown that $\alpha_i = 0$, $i = r_0+1, \ldots, r$ by considering $(\Id - K)u$.
This shows $u \in \Ker K \cap Y \oplus \Range  K \cap Y$.
\end{proof}
\begin{prop}\label{charadmiss}
If the graph $\mG$ is connected, then the following assertions are equivalent.
\begin{enumerate}[(a)]
\item\label{charadmiss:a} The projection $\mathcal {\mathcal P}_K$ is admissible.
\item\label{charadmiss:b} The range of $\It$ is invariant under $\K$, i.\,e., $\K\Range \It \subset \Range \It.$
\item\label{charadmiss:c} There exists a basis of $\Range \It$ consisting of eigenvectors of $\K$.
\end{enumerate}
\end{prop}
\begin{proof}
We start by proving the equivalence of (\ref{charadmiss:a}) and (\ref{charadmiss:b}). Recall that for every $f \in V$ there exists a vector $d^f \in \mathbb C^n$ such that
$$(\I^+)^\top d^f=f(0), \qquad (\I^-)^\top d^f=f(1).$$
The admissibility of the projection is equivalent to the fact that for every $f \in V$ there exists a vector $d^{Pf} \in \mathbb C^n$ such that
$$(\I^+)^\top d^{\mathcal {\mathcal P}_Kf}=\mathcal {\mathcal P}_Kf(0)=Kf(0), \qquad (\I^-)^\top d^{\mathcal {\mathcal P}_Kf}=\mathcal {\mathcal P}_Kf(1)=Kf(1).$$
Inserting the first equation into the second and observing that for all $u \in \mathbb C^n$ there exists a function $f \in (H^1(0,1))^m$ which
is continuous in the nodes such that $d^f=u$ one obtains that (\ref{charadmiss:a}) is equivalent to the fact that for all $u\in \mathbb C^n$ there exists $v\in \mathbb C^n$ such that
$$(\I^+)^\top v = K (\I^+)^\top u, \qquad  (\I^-)^\top v = K (\I^-)^\top u,$$
which can equivalently be stated as
$$
\K\Range \It \subset \Range \It.
$$
The first equivalence is now proved. To see the second equivalence, observe first that the existence of the claimed basis is equivalent to $\Range \It$ being decomposable into $\Range \It=(\Ker \K \cap \Range \It ) \oplus (\Range \K \cap \Range \It).$ Now one can apply the Lemma \ref{decomp} setting $Y:=\Range \It$ and $K:=\K.$
\end{proof}

\begin{lemma}\label{disjointdec}
Consider a decomposition $\mG=\mG_1 \cup \mG_2$ into subgraphs such that every node is contained either in $\mG_1$ or
$\mG_2$. On $\mG_1$, fix a non-admissible orthogonal projection $\mathcal P_{K_1}$.
Then the projection $\mathcal P_K$ on $\mG$ defined by
\[
K:=\begin{pmatrix}
K_1&0\\0&\Id
\end{pmatrix}
\]
is not admissible.
\end{lemma}
\begin{proof}
Since $\mathcal P_{K_1}$ is not admissible, there exists a function $f \in V_1$ such that $\mathcal P_{K_1}f \not\in V_1$, i.\,e., such that the continuity condition is violated in a node $\mv_{k_0}$. It is possible to extend the function $f$ to a function $\tilde{f}$ on the whole graph, such that $d^{\tilde{f}}=0$ in all nodes of $\mG_2$. Then the function $\mathcal P_K\tilde{f}$ does not satisfy the continuity condition in $\mv_{k_0}$, either.
\end{proof}

\subsection{Orthogonality condition --- the matrix $C$}
We are now going to characterize the coefficient matrices $C$ which satify the orthogonality condition; in fact, we will show that the orthogonality condition is equivalent to the invariance of the range of $\mathcal {\mathcal P}_K$ under the coefficient matrix $C$.
\begin{prop}\label{formorth}
Let the sesquilinear form $a$ on $X^2$ be defined as in \eqref{form}, with $M=0$.  Then the following assertions are equivalent.
\begin{enumerate}[(a)]
\item The matrix $C$ satisfies the orthogonality condition~\eqref{Cortho} with respect to $\mathcal {\mathcal P}_K$.
\item The range of $K$ is invariant under the action of $C(x)$ for all $x$, i.\,e.,
\begin{equation}\label{invcoeff}
C(x) \Range K \subset \Range K\qquad \mbox{for all } x\in [0,1].
\end{equation}
\end{enumerate}
\end{prop}
\begin{proof}
Since the space $X^2$ can be decomposed into $X^2=\Range \mathcal {\mathcal P}_K \oplus \Range (\Id-\mathcal {\mathcal P}_K) $, the orthogonality condition \eqref{Cortho} is equivalent to $a(\mathcal P_Ku,(\Id-\mathcal P_K) v) = 0$ for all $u,v\in V$. Using the linearity of the derivative and the self-adjointness of the orthogonal projection $K$, one can compute
\begin{eqnarray*}
a(\mathcal P_Ku,(\Id-\mathcal P_K) v)
	& = & \int_0^1\left(C(x)Ku'(x)\mid(\Id-K)v'(x)\right)dx \\
	& = & \int_0^1\left((\Id-K)C(x)Ku'(x)\mid v'(x)\right)dx,
\end{eqnarray*}
where the inner product is the standard inner product in $\mathbb C^m$.
By a localization argument $\int_0^1 ((\Id-K)C(x)Ku'(x)\mid v'(x))dx=0$ holds for every $u, v \in V$ if and only if $(\Id-K)C(x)K=0$ for all $x\in[0,1]$, i.\,e., $C(x)K=KC(x)K$ for all $x\in[0,1]$.
Since $K$ is a projection this is equivalent to condition~\eqref{invcoeff}.
\end{proof}

\subsection{Orthogonality condition --- the matrix $M$}
Next we characterize the orthogonality condition for the matrix $M$, i.\,e., we want to find equivalent conditions to~\eqref{Mortho1}, where \eqref{Mortho1} can alternatively be stated as
\begin{equation}\label{Mortho}
(Md^{\mathcal {\mathcal P}_Kf} \mid d^{(\Id - \mathcal {\mathcal P}_K) g} )=0 \qquad \hbox{for all }f, g \in V.
\end{equation}
If it is satisfied, we will say that \emph{the matrix $M$ satisfies the orthogonality condition with respect to $K$}, frequently omitting any reference to $K$. For these investigations we introduce the matrix
\[
\mathcal M:= \It D^{-1}MD^{-1}\It^\top,
\]
where we denote $D$ the diagonal matrix with diagonal entries $|\Gamma(\mv_k)|$, the degrees of the nodes $\mv_k$.
Please note that the matrix $\mathcal M$ only depends on $M$ and on the graph structure, but does not depend on the orthogonal projection $K$.

\begin{lemma}\label{fnach}
If the graph $\mG$ has no isolated nodes, then the following assertions hold.
\begin{enumerate}
\item\label{fnach:a} $d^f=D^{-1}\It^\top (f(0),f(1))^\top$ for every $f \in V$.
\item\label{fnach:b} $\Range \It=\left\{(f(0),f(1))^\top\in{\mathbb C}^{2m}: f \in V\right\}$.
\end{enumerate}
\end{lemma}

\begin{proof}
First we will prove the formula
\begin{equation}\label{formuladegree+}
 \mathcal I^+ \left(\mathcal I^+\right)^\top= \diag (\Gamma^+(\mv_k))_{k=1,\ldots,n}.
\end{equation}
In fact, 
\[
(\mathcal I^+ \left(\mathcal I^+\right)^\top)_{lk}=\sum_{i=1}^m \mathcal I^+_{li}\mathcal I^+_{ki}.
\]
Since each edge originates from exactly one node, we obtain that $\mathcal I^+_{li}\mathcal I^+_{ki}=0$ for all $ k\neq l$. Thus,
\[
\sum_{i=1}^m \mathcal I^+_{li} \mathcal I^+_{ki}= \begin{cases} \sum_{i=1}^m \left( \mathcal I^+_{ki}\right)^2,\quad &\hbox{if } k=l,\\0,& \mbox{otherwise.} \end{cases}
\]
Since $\mathcal I^+_{ki}$ equals $1$ exactly $\Gamma^+(\mv_k)$ times and equals $0$ otherwise, the proof of formula~\eqref{formuladegree+} is complete. The analogous formula $ \mathcal I^- \left(\mathcal I^-\right)^\top= \diag (\Gamma^-(\mv_k))_{k=1,\ldots,n}$ 
can be proved likewise. As a consequence, we obtain \[D=\mathcal I^+ \left(\mathcal I^+\right)^\top+\mathcal I^- \left(\mathcal I^-\right)^\top= \It^\top\It.\]

To prove \eqref{fnach:a}, let $f \in V$. By definition, there exists $d^f \in \mathbb C^n$ such that
\begin{equation}\label{existd}
\It d^f= \begin{pmatrix} f(0)\\f(1) \end{pmatrix}.
\end{equation}
We show that the vector $D^{-1}\It^\top\begin{pmatrix}
f(0)\\f(1)
                  \end{pmatrix}$ satisfies the condition \eqref{existd} as well. A direct computation shows that
\[
\It D^{-1}\It^\top\begin{pmatrix} f(0)\\f(1) \end{pmatrix}=\It D^{-1}\It^\top \It d^f 
=\It D^{-1}D d^f
=\It d^f
=\begin{pmatrix}
f(0)\\f(1)\end{pmatrix}.
\]
By the uniqueness of $d^f$, the proof is complete.

For \eqref{fnach:b} notice that since $\{d^f: f \in V\}=\mathbb C^n$,
\[
\Range \It= \left\{\It v: v\in \mathbb C^n\right\}=\left\{\It d^f: f\in V\right\}=  \left\{\begin{pmatrix}
f(0)\\f(1)\end{pmatrix}: f \in V\right\}.
\]
This completes the proof.
\end{proof}

\begin{prop}\label{charM}
Assume the orthogonal projection $\mathcal {\mathcal P}_K$ to be admissible. Then the matrix $M$ satisfies the orthogonality condition~\eqref{Mortho} if and only if
\begin{equation}\label{Mchar}
 \Range \mathcal M \K \It \subset \Range \K.
\end{equation}
\end{prop}
\begin{proof}
We will use the orthogonality condition as stated in \eqref{Mortho}. By Lemma~\ref{fnach}.(1), one obtains for all $f, g \in V$
\[
(Md^{\mathcal {\mathcal P}_Kf} \mid d^{(\Id - \mathcal {\mathcal P}_K) g} )= \left(MD^{-1}\It^\top \K\begin{pmatrix}f(0)\\f(1)\end{pmatrix} \mid D^{-1}\It^\top(\Id-\K)\begin{pmatrix}g(0)\\g(1) \end{pmatrix}\right).
\]
By Lemma~\ref{fnach}.(2), the orthogonality condition is then equivalent to
\begin{equation}\label{boh}
 (MD^{-1}\It^\top \K v \mid D^{-1}\It^\top(\Id-\K)w)=0, \quad \forall v,w \in \Range \It.
\end{equation}
Since $\It$ has real entries and $D^{-1}$ and $(\Id-K)$ are self-adjoint,~\eqref{boh} is equivalent to
\[
 ((\Id-\K) \It D^{-1}MD^{-1}\It^\top \K v \mid w)=0, \quad \forall v,w \in \Range \It,
\]
i.\,e., using the definition of $\mathcal M$
\[
 ((\Id-\K) \mathcal M \K v \mid w)=0, \quad \forall v,w \in \Range \It.
\]
This can be equivalently expressed as
\[
 ((\Id-\K) \mathcal M \K P_{\Range \It}v \mid P_{\Range \It}w)=0, \quad \forall v,w \in \mathbb C^m,
\]
where $P_{\Range \It }$ is the orthogonal projection onto the range of $\It$. Since $P_{\Range \It}$ is self-adjoint,
\[
((\Id-\K) \mathcal M \K P_{\Range \It}v \mid P_{\Range \It}v)=
(P_{\Range \It} (\Id-\K) \mathcal M \K P_{\Range \It}v \mid w)
\]
for all $v,w \in \mathbb C^m$. In fact, we have just proved that the orthogonality condition \eqref{Mortho} is equivalent to
\[
P_{\Range \It} (\Id-\K) \mathcal M \K P_{\Range \It} =0,
\]
which is, finally, the same as
\begin{equation}\label{Mequality}
P_{\Range \It} \mathcal M \K P_{\Range \It} =P_{\Range \It}  \K \mathcal M \K P_{\Range \It}.
\end{equation}
The equation \eqref{Mequality} is the key to prove the claim. Because of the admissibility of $K$, $ \Range \K \It \subset \Range \It$ by Proposition \ref{charadmiss}. Moreover, one sees by the definition of $\mathcal M$ that $\Range \mathcal M \subset \Range \It$, which implies $ \Range\K \mathcal M \subset \Range \It$. Considering both inclusions, one obtains that \eqref{Mequality} is equivalent to $\mathcal M \K P_{\Range \It} =  \K \mathcal M \K P_{\Range \It}$. 
In fact, we are asking that $\K$ acts as the identity matrix on $\Range \mathcal M \K P_{\Range \It}=  \Range\mathcal M \K \It$, i.\,e., the orthogonality condition \eqref{Mortho} is equivalent to
$ \Range \mathcal M \K \It \subset \Range \K$. This  concludes the  proof.
\end{proof}

Although it is easy to check those range inclusions numerically for concrete examples, the following sufficient conditions may
be more convenient in some cases.
\begin{cor}
Consider an admissible projection $\mathcal {\mathcal P}_K$. If 
$\Range\mathcal M  \It \subset \Range \K$ or $\Range\mathcal M \K \subset \Range \K$, then $M$ satifies the orthogonality condition.
\end{cor}
\begin{proof}
To see that $\Range\mathcal M  \It \subset \Range \K$ is sufficient, observe that admissibility of $\mathcal P_K$ implies $\Range\K \It \subset \Range \It$,
and hence condition \eqref{Mchar} is fulfilled.
Moreover, $\Range\mathcal M \K \subset \Range \K$ is also sufficient by a similar argument, since $ \Range\K \It \subset \Range \K$.
\end{proof}

As an application we give a a simpler characterization for a special subspace of $\mathbb C^m$ in the case of a bipartite graph.
More precisely, we consider the smallest subspace of $\mathbb C^m$ whose orthogonal projection $K$ satisfies $K\1 = \1$, i.\,e.,
$Y = \left\{ (c, c, \ldots, c)^T \mid c \in \mathbb C \right\}$.
\begin{prop}\label{Morthobip}
Let the graph $\mG$ be bipartite. Then $M$ satisfies the orthogonality condition with respect to the orthogonal projection $K:=\left(\frac{1}{m}\right)_{i,j=1,\ldots,m}$ if and only if there exist values $(\alpha_{ij})_{i,j=1,2}$ such that
\[ \begin{array}{lll}
\alpha_{11}{|\Gamma(\mv_\ell)|}=\sum_{k=1}^{n_1} m_{\ell k}, & \alpha_{12}{|\Gamma(\mv_\ell)|}=\sum_{k=n_1+1}^{n} m_{\ell k } &  \hbox{for all }\ell=1,\ldots,n_1, \quad \text{and} \\
\alpha_{21}{|\Gamma(\mv_\ell)|}=\sum_{k=1}^{n_1}m_{\ell k}, & \alpha_{22}{|\Gamma(\mv_\ell)|}=\sum_{k=n_1+1}^{n} m_{  \ell k} & \hbox{for all }\ell=n_1+1,\ldots,n.
\end{array} \]
\end{prop}

\begin{proof}
First we show $\Range \K \It = \Range K$, and then we prove a characterization of those matrices $\mathcal M$ that leave invariant the subspace $\Range \K$. Then the claim will follow by Proposition~\ref{charM}.

Without loss of generality we may assume that for $k=1,\ldots,n_1$ the node $\mv_k$ has only outgoing edges, i.\,e., $\Gamma(\mv_k) = \Gamma^+(\mv_k)$,
and that for $k=n_1+1,\ldots,n$ the node $\mv_k$ has only incoming edges, i.\,e., $\Gamma(\mv_k) = \Gamma^-(\mv_k)$,
reordering the nodes otherwise. We are going to prove
\begin{equation}\label{bipKI}
\Range \K \It=\Range \K=\left< \mathbb 1_{\{1,\ldots,m\}},\mathbb 1_{\{m+1,\ldots,2m\}} \right>.
\end{equation}
The second equality follows from the definition of $K$. Moreover, $\Range \K \It \subset \Range \K$ is obvious.
Since there exists $f \in V$ such that $d^f = \1_{\{1, \ldots, n_1\}}$, 
Lemma~\ref{fnach}.\,(2) implies $\mathbb 1_{\{1,\ldots,m\}} \in \Range \It$.
Analoguously $d^f = \1_{\{n_1+1, \ldots, n\}}$ yields $\mathbb 1_{\{m+1,\ldots,2m\}} \in \Range \It$.
We have already observed that $\mathbb 1_{\{1,\ldots,m\}}$ and $\mathbb 1_{\{m+1,\ldots,m\}}$ are invariant under $\K$, which implies
$$\left<\mathbb 1_{\{1,\ldots,m\}},\mathbb 1_{\{m+1,\ldots,m\}} \right> \subset \Range \K \It,$$
and this proves the claim \eqref{bipKI}.

Next we characterize the matrices $\mathcal M$ which leave $\Range \K$ invariant. For this we use the bipartite decomposition of
$$M^w:=D^{-1}MD^{-1}= \left(\frac{m_{ij}}{|\Gamma(\mv_k)||\Gamma(\mv_\ell)|}\right)_{k,\ell=1,\ldots,n}$$
which is induced by the bipartite decomposition of the graph $\mG$, i.\,e., we write
$$
M^w=\begin{pmatrix}
M^w_{11} & M^w_{12}\\
M^w_{21} & M^w_{22}
    \end{pmatrix},
$$
where $M^w_{11} \in M_{n_1, n_1}$, $M^w_{12} \in M_{n_1, n-n_1}$, $M^w_{21} \in M_{n-n_1, n_1}$, and $M^w_{22} \in M_{n-n_1, n-n_1}$.

Moreover, since $\mG$ is bipartite, the incidence matrices decompose into
$$
\mathcal I^{+}=\begin{pmatrix} \mathcal I^{+}_1 \\ 0 \end{pmatrix} \mbox{ and }
\mathcal I^{-}=\begin{pmatrix} 0      \\ \mathcal I^{-}_2 \end{pmatrix},
$$
where $\mathcal I_1^{+} \in M_{n_1, m}$ and $\mathcal I_2^{+} \in M_{n-n_1, m}$.

We will use the above decompositions in order to obtain a useful formula for $\mathcal M$. Using the definition of $\It$, as formulated in \eqref{deftilde}, we first compute
$$
\mathcal M= \It M^w \It^\top =\It(M^w \mathcal I^+,M^w \mathcal I^-)=
\begin{pmatrix}
\mathcal I^{+^\top}M^w \mathcal I^+  &\mathcal I^{+^\top}M^w \mathcal I^-\\
\mathcal I^{-^\top}M^w \mathcal I^+  &\mathcal I^{-^\top}M^w \mathcal I^-
\end{pmatrix},
$$
which yields, inserting the decompositions,
$$\mathcal M=
\begin{pmatrix}
\mathcal I_1^{+^\top}M^w_{11} \mathcal I^+_1  &\mathcal I_1^{+^\top}M^w_{12} \mathcal I^-_2\\
\mathcal I_2^{-^\top}M_{21}^w \mathcal I_1^+  &\mathcal I_2^{-^\top}M_{22}^w \mathcal I_2^-
\end{pmatrix}.
 $$
It is possible to identify the block-matrices appearing in the above expression. In fact, the following identities hold.
\begin{eqnarray*}
\mathcal I_1^{+^\top}M^w_{11} \mathcal I^+_1 & = & \left( \frac{m_{\me_i(0)\me_j(0)}}{|\Gamma(\me_i(0)| |\Gamma(\me_j(0))|} \right)_{i,j=1,\ldots,m} \\
\mathcal I_2^{-^\top}M^w_{22} \mathcal I^-_2 & = & \left( \frac{m_{\me_i(1)\me_j(1)}}{|\Gamma(\me_i(1)| |\Gamma(\me_j(1))|} \right)_{i,j=1,\ldots,m} \\
\mathcal I_1^{+^\top}M^w_{12} \mathcal I^-_2 & = & \left( \frac{m_{\me_i(0)\me_j(1)}}{|\Gamma(\me_i(0)| |\Gamma(\me_j(1))|} \right)_{i,j=1,\ldots,m} \\
\mathcal I_2^{-^\top}M^w_{21} \mathcal I^+_1 & = & \left( \frac{m_{\me_i(1)\me_j(0)}}{|\Gamma(\me_i(1)| |\Gamma(\me_j(0))|} \right)_{i,j=1,\ldots,m}
\end{eqnarray*}
Here we write $m_{\me_i(0) \me_j(0)}:=m_{k\ell}$ if $\me_i(0)=\mv_k$, $\me_j(0)=\mv_\ell$, and analogously for $\me_j(1)$.

We have already observed that $\{\mathbb 1_{\{1,\ldots,m\}},\mathbb 1_{\{m+1,\ldots,2m\}}\}$ is a basis of $\Range \K$. Using the above decompositions, one can compute
$$
\mathcal M \mathbb 1_{\{1,\ldots,m\}}=
\begin{pmatrix}
\left(\mathcal I_1^{+^\top}M^w_{11} \mathcal I^+_1 \right)\mathbb 1_{\mathbb C^m}  \\
\left(\mathcal I_2^{-^\top}M^w_{21} \mathcal I^+_1 \right)\mathbb 1_{\mathbb C^m}
\end{pmatrix}
=
\begin{pmatrix}
\left(\sum_{j=1}^m \frac{m_{\me_i(0)\me_j(0)}}{|\Gamma(\me_i(0))| |\Gamma(\me_j(0))|} \right)_{i=1,\ldots,m}\\
\left(\sum_{j=1}^m \frac{m_{\me_i(1)\me_j(0)}}{|\Gamma(\me_i(1))| |\Gamma(\me_j(0))|} \right)_{i=1,\ldots,m}
\end{pmatrix}.
$$
In these sums each $\me_j(0)$ appears exactly $|\Gamma^+(\me_j(0))|$ times. Collecting the same summands we may write
$$
\mathcal M \mathbb 1_{\{1,\ldots,m\}}=
\begin{pmatrix}
\left( \sum_{k=1}^n |\Gamma^+(\mv_k)| \frac{m_{\me_i(0)k}}{|\Gamma(\me_i(0))| |\Gamma(\mv_k)|} \right)_{i=1,\ldots,m}\\
\left( \sum_{k=1}^n |\Gamma^+(\mv_k)| \frac{m_{\me_i(1)k}}{|\Gamma(\me_i(1))| |\Gamma(\mv_k)|} \right)_{i=1,\ldots,m}
\end{pmatrix}.
$$
In fact, since $|\Gamma^+(\mv_k)|$ appears as a factor, only for the first $n_1$ vertices the summand does not vanish. Thus, we see that
$$
\mathcal M \mathbb 1_{\{1,\ldots,m\}}=
\begin{pmatrix}
\left(\sum_{k=1}^{n_1} |\Gamma^-(\mv_k)| \frac{m_{\me_i(0)k}}{|\Gamma(\me_i(0))| |\Gamma(\mv_k)|} \right)_{i=1,\ldots,m}\\
\left(\sum_{k=1}^{n_1} |\Gamma^-(\mv_k)| \frac{m_{\me_i(1)k}}{|\Gamma(\me_i(1))| |\Gamma(\mv_k)|} \right)_{i=1,\ldots,m}
\end{pmatrix}.
$$
Since $|\Gamma^+(\mv_k)|=|\Gamma(\mv_k)|$ for $k=1,\ldots,n_1$,
$$
\mathcal M \mathbb 1_{\{1,\ldots,m\}}=
\begin{pmatrix}
\left(\sum_{k=1}^{n_1}  \frac{m_{\me_i(0)k}}{|\Gamma(\me_i(0))|} \right)_{i=1,\ldots,m}\\
\left(\sum_{k=1}^{n_1}  \frac{m_{\me_i(1)k}}{|\Gamma(\me_i(1))|} \right)_{i=1,\ldots,m}
\end{pmatrix}.
$$
We can easily check whether $\mathcal M \mathbb 1_{\{1,\ldots,m\}} \in \Range \K$. Using \eqref{bipKI} one sees that this is the case if and only if sums above do not depend on $\me_i$,
i.\,e., $\mathcal M \mathbb 1_{\{1,\ldots,m\}} \in \Range \K$ if and only if there exist $\alpha_{11}, \alpha_{21} \in \mathbb C$ such that 
$ \alpha_{11}{|\Gamma(\mv_\ell)|}=\sum_{k=1}^{n_1} m_{\ell k}$ for all $\ell=1,\ldots,n_1$, and
$ \alpha_{21}{|\Gamma(\mv_\ell)|}=\sum_{k=1}^{n_1}m_{\ell k} $ for all $\ell=n_1+1,\ldots$.
By a similar computation, one can also see that $\mathcal M \mathbb 1_{\{m+1,\ldots,2m\}} \in \Range \K$ if and only if there exist
$\alpha_{12}, \alpha_{22} \in \mathbb C$ such that
$ \alpha_{12}{|\Gamma(\mv_\ell)|}=\sum_{k=n_1+1}^{n} m_{\ell k }$  for all $\ell=1,\ldots,n_1$, 
and $ \alpha_{22}{|\Gamma(\mv_\ell)|}=\sum_{k=n_1+1}^{n} m_{  \ell k}$ for all $\ell=n_1+1,\ldots,n$.
This completes the proof.
\end{proof}

\begin{exa}\label{stochasticmat}
Consider a regular, bipartite graph $\mG$, with the bipartite node decomposition $\mG=\{\mv_1,\ldots,\mv_{n_1}\}\cup\{\mv_{n_1+1},\ldots,\mv_{n}\}$. Set $n_2:=n-n_1$ and consider row-stochastic matrices $M_{ij} \in M_{n_i, n_j}$. Then for arbitrary $\alpha_{ij} \in \mathbb C$ all matrices of the form
$$
M:=\begin{pmatrix}
\alpha_{11} M_{11} & \alpha_{12} M_{12}\\
\alpha_{21} M_{21} & \alpha_{22} M_{22}
   \end{pmatrix},
$$
satisfy the orthogonal condition with respect to $K$ defined as in Proposition~\ref{Morthobip}.
\end{exa}

\section{Classes of graphs}
In this section we will discuss some classes of graphs, combining the results of the preceeding sections. We present some (non-standard) graph theoretical definitions we will use through this section.

\begin{defi}\label{graphdef}
Let $\mG$ a graph with no isolated nodes, i.\,e., such that
$\Gamma(\mv_k)\geq 1$ for all $k=1,\ldots,n$.
\begin{itemize}
\item We call the graph $\mG$ \emph{completely unconnected} if $\mG$ is the union of disjoint compact intervals, i.\,e., if $\mG$ is a regular graph of degree $1$.
\item We call the graph $\mG$ an \emph{inbound} (respectively, \emph{outbound}) star, if there exists a node $\mv_k$ such that
$\me_j(1)=\mv_k$, (respectively, if $\me_j(0)=\mv_k$), for all $j=1,\ldots,m$. We call the graph $\mG$ a star if it is an inbound or outbound star and $\mv_k$ the \emph{center}
of the star.
\item We call the graph $\mG$ \emph{bipartite} if each node has only either incoming or outgoing edges.
\item We call the graph $\mG$ \emph{Eulerian} if all nodes have the same number of incoming and outgoing edges.
\item We call a graph $\mG$ a \emph{layer graph} if there  exist disjoint sets $V_1,\ldots,V_L$ such that 
\begin{itemize}
\item $V=\cup_{p=1}^L V_p$,
\item $\me_j(0) \in V_p$ implies $\me_j(1) \in V_{p+1}$  for all $p=1,\ldots,L-1$, and 
\item $\me_j(0) \in V_L$ implies $\me_j(1) \in V_{1}$.
\end{itemize}
 Nodes belonging to $V_p$ are said to lie in the $p$\textsuperscript{th}
layer. Edges outgoing from nodes in the $p$\textsuperscript{th} layer are also said to
lie in the $p$\textsuperscript{th} layer.
\item We call a layer graph \emph{symmetric} if the incoming and outgoing degrees of the nodes only depends on the layer, i.\,e., if there exist numbers $I(p), O(p) \in \mathbb N_0$ such that $|\Gamma^+(\mv)|=I(p), |\Gamma^-(\mv)|=O(p)$ for all nodes $\mv$ in the $p$\textsuperscript{th} layer.
\end{itemize}
\end{defi}

\subsection{Bipartite and Euler Graphs}
It is possible to characterize some classes of graphs by the
admissibility of the matrix from Proposition~\ref{Morthobip}.
\begin{theo}\label{mittelunggraph}
Consider the orthogonal projection $K$ defined by
\begin{equation}\label{averaging}
K:=\left(\frac{1}{m}\right)_{i,j=1,\ldots,m}.
\end{equation}
Then $\mathcal P_K$ is admissible if and only if $\mG$ is bipartite or Eulerian.
\end{theo}

\begin{proof}
Fix $f \in V$ and observe that $\mathcal P_Kf$ always lies
in $(H^1(0,1))^m$ since every component is a linear combination of
$H^1$ functions. So $V$ is invariant if and only if $\mathcal P_kf$ is
continuous in the nodes.
Let $\mV_1 \subset \mV$ denote the set of all vertices
having outgoing edges, and let $\mV_2 \subset \mV$
denote the set of all vertices having incoming edges. We distiguish
two cases. First, assume  $\mV_1 \cap \mV_2 = \emptyset$.
Then $\mG$ is a bipartite graph.

On the other hand, if $\mV_1 \cap \mV_2  \not= \emptyset$, then
by definition of $K$ a vector
$d^{\mathcal {\mathcal P}_Kf}$ exists if and only if
\begin{equation} \label{sameaverage}
\sum_{j=1}^m\frac{f_j(0)}{m} = \sum_{j=1}^m\frac{f_j(1)}{m}.
\end{equation}
We show now that the equality
\eqref{sameaverage} is equivalent to the graph being Eulerian.
First, assume that \eqref{sameaverage} holds for every $f \in V$.
Fix an arbitrary $\mv_k \in \mV$ and choose $f \in V$ such
that $d^f = \1_{\{i\}}$. Then
$$
\frac{1}{m} \left|\Gamma^+(\mv_k)\right| = \sum_{j=1}^m\frac{f_j(0)}{m} =
\sum_{j=1}^m\frac{f_j(1)}{m} = \frac{1}{m} \left|\Gamma^-(\mv_k)\right|.
$$
Thus it is necessary that $\left|\Gamma^-(\mv_k)\right| = \left|\Gamma^+(\mv_k)\right|$ holds
for every $k=1,\ldots,n$. Conversely, assume that
$\left|\Gamma^-(\mv_k)\right| = \left|\Gamma^+(\mv_k)\right|$ holds for every $k=1,\ldots,n$.
Then
$$
\sum_{j=1}^m\frac{f_j(0)}{m} =
\frac{1}{m}\sum_{k=1}^{n}\left|\Gamma^+(\mv_k)\right|d^f_k =
\frac{1}{m}\sum_{k=1}^{n}\left|\Gamma^-(\mv_k)\right|d^f_k =
\sum_{j=1}^m\frac{f_j(1)}{m}.
$$
Hence \eqref{sameaverage} is satisfied, so this condition is also sufficient.

It only remains to show that indeed for every bipartite graph $K$ is admissible.
To see this, note that for an arbitrary $f \in V$ the vector $d^{\mathcal P_Kf}$ can
be chosen to equal $\sum_{i=1}^m\frac{f_i(0)}{m}$ in all components belonging to nodes in
$\mV_1$ and to equal $\sum_{i=1}^m\frac{f_i(1)}{m}$ in all
components belonging to $\mV_2$. This shows continuity of $\mathcal P_Kf$ in the nodes,
thus implying $\mathcal P_Kf \in V$.
\end{proof}

\begin{rem}
The matrix $K$ defined in \eqref{averaging} acts on a vector $v\in
\mathbb C^m$ by substituting each component by the average of all
components of the vector. The range of such a matrix is thus
one-dimensional, and one sees that
$$
\Range \mathcal {\mathcal P}_K=\{f \in V: f_i=f_j \mbox{ for all }
i,j=1,\ldots,m\}.
$$
Such functions are symmetric in the sense that they are equal on each
edge at the same point of the parametrization. In fact, Theorem \ref{mittelunggraph} characterizes the
admissibility of projections whose ranges consist of the functions that are
symmetric on the network. It thus gives a first answer to the problem
stated in Remark \ref{motivexa}.
\end{rem}

\subsection{Stars}
Main result of this subsection is a characterization of stars in the class of the simple graphs.
We first investigate the admissibility of projections.

\begin{prop}\label{unconnstar}
The following assertions hold.
\begin{enumerate}
\item The graph $\mG$ is completely unconnected if and only if
$\mathcal {\mathcal P}_K$ is admissibile for all orthogonal
projections $K$.
\item Let $\mG$ be a simple, connected graph. Then $\mG$
is a star if and only if $\mathcal {\mathcal P}_K$ is admissibile for
all orthogonal projections $K$ with eigenvector $\1$.
\end{enumerate}
\end{prop}

\begin{proof}
(1) Since the graph $\mG$ is completely unconnected, the
continuity condition in $V$ is empty, and therefore each $\mathcal
{\mathcal P}_K$ is admissible. Conversely, if $\mG$ is not
completely unconnected, then it is possible to decompose $\mG$
into the disjoint union of a connected graph $\mG_1$ with $m_1$
edges and the remaining graph $\mG_2$. Let $K_1$ be an
orthogonal projection of $\mathbb C^{m_1}$, which does not have $\1$ as
an eigenvector. Lemma~\ref{1eigen} and Lemma~\ref{disjointdec} assert
that the orthogonal projection
$$\begin{pmatrix} K_1&0\\0&\Id \end{pmatrix}$$ 
is not admissible.

(2) Without loss of generality, we prove the claim for an outgoing
star with center $\mv_1$ and with the natural numbering of the other
nodes. Let the graph $\mG$ be a star and $K$ be a projection
such that $K\1=\1$. In fact, for this star
$$\It=\begin{pmatrix}\1&0\\0&\Id_m\end{pmatrix}.$$ 
Since now $K$ has
$\1$ as eigenvector to the eigenvalue 1, one can compute
$$
\K\It=\begin{pmatrix}  K\1 &0 \\0 & K \Id_m\end{pmatrix} =
\begin{pmatrix}  \1 &0 \\0 & K \end{pmatrix}.
$$
It is now clear that $\Range\K\It \subset \Range\It $, and this
implies the admissibility of $\mathcal {\mathcal P}_K$. Conversely,
assume that the graph $\mG$ is not a star. One sees that this implies
the existence of an undirected path of length
$3$. We will denote it by
$\me_1,\me_2,\me_3$, possibly relabelling the edges.
Our strategy is the following: for each graph that is a path consisting of $3$ edges we construct a non-admissible projection $\mathcal P_L$ where $L\1 = \1$.
We then consider the projection $\mathcal P_K$, where $K$ is
\[
K:=\begin{pmatrix}
L & 0\\0&\Id
  \end{pmatrix}.
\]
Then, by Lemma \ref{disjointdec}, we conclude that $\mathcal P_K$ is not admissible,
although $\1$ s an eigenvector of $K$.

First, consider cycles of length $3$. Since each edge can be
directed arbitrarily, there are $8$ such graphs. Let us
start with the case of a not strongly connected graph.
Such graphs are neither Eulerian nor bipartite.
Thus, Theorem~\ref{mittelunggraph} provides an example of an $L$ as requested.
If the graph is a (directed) cycle such that
$\me_1(0)=\mv_1$, consider the projection
\[ L:=\begin{pmatrix}
\frac{1}{2} & \frac{1}{2} & 0\\
\frac{1}{2} & \frac{1}{2} & 0\\
  0&0&1
\end{pmatrix}\]
and the function $f$ defined by $f(x):=(x,1-x,0)^\top \in V$. One sees that $f \in V$ but
$\mathcal P_Kf \not\in V$, since $\mathcal P_Kf(x) = (\frac{1}{2}, \frac{1}{2},0)^\top$ for a.\,e. $x\in(0,1)$.

Consider now the lines of length $3$. We split this into three possible
cases: $\mG$ may be bipartite line,
a (directed) line,
or neither a (directed) line nor a bipartite graph.
In the last two cases the graphs is neither bipartite nor
Eulerian, and hence we can use Theorem~\ref{mittelunggraph} again.
In the case of a bipartite line, let us consider the projection
\[
L:=\begin{pmatrix}
\frac{1}{2} & \frac{1}{2} & 0\\
0&0&1\\
\frac{1}{2} & \frac{1}{2} & 0\\
\end{pmatrix}
\]
for the parametrization
$\me_1(0)=\mv_1$, $\me_1(1)=\me_2(1)=\mv_2$, $\me_2(0)=\me_3(0)=\mv_3$, and
$\me_3(1)=\mv_4$. Consider the function $f(x):=(x,x,0)^\top$. Again, $f \in V$ but  $\mathcal P_Kf\not\in V$, since $\mathcal P_Kf(x)= (\frac{x}{2}, x,\frac{x}{2})^\top$ for a.\,e. $x\in(0,1)$. This completes the proof.
\end{proof}

\begin{rem}
In Proposition \ref{unconnstar}, (2) we have assumed the graph $\mG$ to have no multiple edges. In fact, it is not possible to relax this
condition, since all orthogonal projections with eigenvector $\mathbb 1$ are admissible on all connected graphs consisting of $2$ nodes and
 $m$ edges for each $m \in \mathbb N$ and each orientation of the
edges.
\end{rem}

Now we investigate the orthogonality condition for diagonal matrices $C$.
This will show that for a wide class of matrices $C$
there cannot exist non-trivial invariant subspaces of the form considered in this paper.

\begin{lemma}\label{diagonalmatrix}
Let ${\mathcal D}$ be a constant diagonal matrix with entries $d_i>0$. Then the following assertions hold.
       \begin{enumerate}
       \item
			Assume the coefficients $d_i$ to be pairwise different. If $K$ is
			an orthogonal projection with eigenvector $\1$ such that ${\mathcal D} \Range K
			\subset \Range K$, then $K$ is trivial, i.\,e., $K = \Id$ or $K = 0$.
       \item
               Assume that there exists $i_0 \neq j_0$ such that $d_{i_0} = d_{j_0}$.
               Then there exists a nontrivial orthogonal projection $K$ with eigenvector $\1$
               such that ${\mathcal D} \Range K \subset \Range K$.
       \end{enumerate}
\end{lemma}

\begin{proof}   Observe that it is possible to compute the powers of ${\mathcal D}$
explicitly, since it is diagonal. In fact,
               ${\mathcal D}^k=\diag(c^k_i)_{i=\ldots,m}$ for every $k \in \mathbb N_0$. 

(1) Since $K$ is an orthogonal projection and $\1$ is an eigenvector, either $K\1=\1$ or $K\1=0$. If $K \1 = \1$, i.\,e., $\1 \in \Range K$,  we see by induction that
$(d_1^k, d_2^k, \dots, d_m^k) = {\mathcal D}^k \1 \in \Range K$
               for every $k \in \mathbb{N}$ since $\Range K$ is invariant
under the action of ${\mathcal D}$.
               Now, the matrix $V := (d_{ij})_{i,j=1, \dots, m}$, defined by
               $$d_{ij}:=d_i^{j-1} \qquad {i,j=1, \dots, m}$$ is the Vandermonde
matrix induced by the vector $(d_i)_{i=1,\ldots,m}$, which is regular
since the $d_i$ are pairwise different. From this we see $\Range K =
\mathbb{C}^k$, i.\,e., $K = \Id$.

If on the other hand $K \1 = 0$, then fix $v \in \Range K$. Since $\1$ is in the
kernel of $K$, $\Range K\subset \left<\1\right>^\bot$.
Since the range of $K$ is invariant under the
action of the matrix $C$, we obtain
$\left( C^nv \mid \1 \right) = \sum_{i=1}^m d_i^n v_i = 0$ for every $n \in \mathbb{N}$. In
particular, $v$ satisfies the equation $V^T v = 0$. Since $V$ is
regular, we obtain $v=0$, which implies $\Range K = \{0\}$, hence $K =
0$.

(2) In order to prove the second assertion, let $i \neq j$ such that $d_i = d_j$.
               Consider $$Y := \Span\{ {\mathcal D}^n \1 = (d_1^n, \dots, d_m^n) \mid n \in
\mathbb{N} \} \subset \mathbb{C}^m,$$ and let $K$ be the orthogonal
projection onto $Y$. Since $\1 \in Y$, $K\1=\1$, and in
particular $K\not=0$.
               Moreover, $v_i=v_j$ for all $v \in Y$, which implies
$\Range K \not=\mathbb C^m$. As a consequence $K\not=\Id$. Finally,
${\mathcal D}Y \subset Y$, and hence the range of $K$ is invariant under
the action of ${\mathcal D}$. This completes the proof.
\end{proof}

Combining the previous two statements we deduce the following.
\begin{prop}\label{star}
Let the graph $\mG$ be connected. If the coefficient matrix $C$ is diagonal and  $M = 0$, then the following assertions hold.
\begin{enumerate}
\item  Let $C$ be constant and $\mG$ be a star. If there exist $i_0,j_0$ such that $c_{i_0}=c_{j_0}$, then the subspace 
$$\mathcal Y:=\{ f \in X^2 \mid f_{i_0}(x)=f_{j_0}(x)  \mbox{\ for
a.\,e. } x \in (0,1) \}$$
is invariant under the action of $(e^{tA})_{t\geq 0}$.
\item   Assume the coefficients $c_i(x_0)$ to be pairwise different for some $x_0\in[0,1]$. If $Y$ is a nontrivial linear subspace of $\mathbb C^m$, then $\mathcal Y$,
defined as above, is not invariant under the action of $(e^{tA})_{t\geq 0}$.
\end{enumerate}
\end{prop}

\begin{proof}
       (1) Without loss of generality, assume $c_1=c_2$. Consider the
subspace $$Y:=\{v\in \mathbb C^m: v_1=v_2\}$$ and let $K$ be the
orthogonal projection onto $Y$. Since $\1 \in Y$, $K\1=\1$.
Furthermore, by Proposition \ref{unconnstar} $K$ is admissible, since
$\mG$ is a star. Let $v \in Y$. Computing now
$Cv=(c_jv_j)_{j=1,\ldots,m}$ shows that $Cv \in Y$, since $c_1=c_2$
and $v_1=v_2$. This shows that $C\Range K \subset \Range K$, thus
completing the proof of the first claim.

(2) Let $Y \subset \mathbb{C}^m$ be a linear subspace, and let $K$ be the orthogonal projection onto $Y$.
       Remember that the invariance of this subspace is equivalent to the fact that
       $K$ is admissible and that the sesquilinear form $a$ satisfies the
orthogonality condition with respect to $\mathcal {\mathcal P}_K$. If $K$ is not admissible, then the proof is
complete. Thus, assume that $K$ is admissible. Since the graph $\mG$ is assumed to be connected, $\1$ is an eigenvector of $K$.
The orthogonality condition is equivalent to $C(x)Y \subset Y$ for all $x\in[0,1]$, according to Proposition \ref{formorth}. Using Lemma \ref{diagonalmatrix} we see that since the diagonal entries of $C(x_0)$ are
pairwise different, this is not possible for non-trivial $K$.
Hence the proof is complete.
\end{proof}

\subsection{Layer Graphs}
In this section we prove an admissibility result for symmetric layer graphs.
We start fixing a canonical numbering of the edges of a layer graph.
First observe that the node decomposition induces an edge
decomposition $E=\cup_{p=1}^L E_p$ by setting
$$
E_p:=\{\me \in\mE : \me \mbox{ lies in the $p$\textsuperscript{th} layer}\}.
$$
After relabeling the edges we may assume that there exist $L_p$,
$p=1,\ldots,L+1$ satisfying
\begin{enumerate}
\item $L_1=0$;
\item $\me_i(0)=\me_j(0) \mbox{ or } \me_i(1)=\me_j(1)$ implies $L_{p-1} < i,j \leq L_{p}$ for some $p$;

\item $\me_i(0)=\me_j(1)$ implies $L_{p-1} < j \leq L_p < i \leq
L_{p+1}$ for some $p$.
\end{enumerate}
The numbering obtained in such a way has the property that $\me_i$ is
in the $p$\textsuperscript{th} layer if and only if $L_p < i \leq L_{p+1}$. In fact,
all edges $\me_i$ such that $i \leq L_{p+1}$ are in any of the first
$p$ layers.

We are going to exhibit a class of admissible projections. Altough the
result is not a complete characterization, it is optimal in a
sense we will explain later.
\begin{prop}\label{layer}
Consider a symmetric layer graph $\mG$ and the orthogonal projection $K$
\begin{equation}\label{layerproj}
K=\begin{pmatrix}
(\frac{1}{|E_1|})_{i,j=1,\ldots,|E_1|} &&0\\
&\ddots&\\
0&&(\frac{1}{|E_L|})_{i,j=1,\ldots,|E_L|}
  \end{pmatrix},
\end{equation}
where $|E_p|$, $p=1,\ldots,L$ denotes the number of edges in the
$p^{\rm th}$ layer. Then $\mathcal P_K$ is admissible.
\end{prop}

\begin{proof}
One has to check the continuity condition for each $p=1,\ldots, L-1$ in every node of the $p$\textsuperscript{th} layer. Define the auxiliary
function
$$
\lambda: k \mapsto \mbox{layer of the node $\mv_k$}.
$$
We thus have to check continuity in those nodes $\mv_k$ such that $\lambda(k) = p$, $p=1,\ldots, L-1$.

The set $\lambda^{-1}(p)$ can be represented in the form
$$
\lambda^{-1}(p)=\{k: \exists i \in \{L_p+1,\ldots,L_{p+1}\} \mbox{ s.\,t.
}\me_i(1)=\mv_k\},
$$
as well as in the form
$$
\lambda^{-1}(p)=\{k: \exists i \in \{L_{p+1}+1,\ldots,L_{p+2}\} \mbox{
s.\,t. }\me_i(0)=\mv_k\},
$$
whenever the expression is defined.
By the definition of $K$, one sees that for all $p=1,\ldots,L-1$ and
all $i,j=L_p+1,\ldots,L_{p+1}$ the identities
\begin{equation}\label{indimitt1}
\mathcal {\mathcal P}_Kf_i(1)=\mathcal {\mathcal P}_Kf_j(1), \qquad
\mathcal {\mathcal P}_Kf_i(0)=\mathcal {\mathcal P}_Kf_j(0).
\end{equation}
hold.
As a consequence, for layers having incoming or outgoing degree $0$,
the continuity is obvious. Assume now that $I(p)\neq 0$ and
$O(p)\neq 0$.

For an edge $\me_i$ in the $p$\textsuperscript{th} layer and for $f \in \mathcal V$,
\begin{eqnarray*}
\mathcal ({\mathcal P}_K f)_i(1)=\sum_{i=L_p+1}^{L_{p+1}} \frac{ f_i(1)}{|E_p|}
=\sum_{k\in \lambda^{-1}(p)} |\Gamma^-(\mv_k)|\frac{ f_i(1)}{|E_p|}.
\end{eqnarray*}
Recall that since our graph is symmetric, the incidence degree
$|\Gamma^-(\mv_k)|$ only depends on the layer, and therefore we can write
$$
\mathcal ({\mathcal P}_K f)_i(1)=\sum_{k\in \lambda^{-1}(p)} |I(p)|\frac{
f(\mv_k)}{|E_p|}
=\frac{|I(p)|}{|E_p|}\sum_{k\in \lambda^{-1}(p)} f(\mv_k).
$$
With analogous computations we obtain for edges in the $p+1$ layer
$$
\mathcal ({\mathcal P}_K f)_i(0)=\frac{|O(p)|}{|E_{p+1}|}\sum_{k\in
\lambda^{-1}(p)} f(\mv_k).
$$
Observe that the identities
$|E_{p+1}|=|\lambda^{-1}(p)||O(p)|$ and $|E_p|=|\lambda^{-1}(p)||I(p)|$
imply
${|O(p)|}{|E_{p+1}|}^{-1}={|I(p)|}{|E_p|}^{-1}$.
We have thus proved that $
(\mathcal P_K f)_i(1)=(\mathcal P_K f)_j(0)$ for all $\me_i,\me_j$ such that $i \in
\lambda^{-1}(p), j\in \lambda^{-1}(p+1)$. This completes the proof.
\end{proof}

\begin{cor}
Let $\mG$ be a symmetric layer graph. If $M=0$ and $C=c(x)\Id$ for some function $0<c\in C^1[0,1]$, then the space
$$
\mathcal Y:=\{f \in L^2: f_i=f_j \mbox{ for all } i,j \in \ell^{-1}(p), \; p=1,\ldots,L\}
$$
is invariant under the action of $(e^{tA})_{t\geq 0}$.
\end{cor}

\begin{rems}\mbox{}
\begin{enumerate}[(1)]
\item
	The class of the \emph{layer graphs} is not a common object in the graph theoretical literature. 
	In fact, layer graphs are nothing but (directed) $p$-partite
	graphs, for which collapsing the components of the graphs to
	a single vertex leads to a finite line or to a cycle.  In particular, homogeneous trees of finite depth are symmetric
	layer graphs. Such graphs play a role in the investigation of
	biological neural networks.
\item
	The symmetry condition in Proposition~\ref{layer} cannot be
	relaxed. To see this, consider the following simple example. Let $\mG$
	be an outgoing star of order two and consider two copies of $\mG$.
	Identifying two of the external nodes defines a layer graph. One can
	show that the orthogonal projection defined in~\eqref{layerproj} is not admissble, due to the two free nodes in the second layer.
\item
	It seems to be possible to extend the result of Proposition~\ref{layer} to non-symmetric layer graphs, requiring some weaker condition and suitably weighting the projection of~\eqref{layerproj} according to the degrees. However, such results are quite techical. Presenting them in detail goes beyond the scope of this paper.
\end{enumerate}
\end{rems}

\section{Applications}

\subsection{Ephaptic coupling of biological fibers}

In the modern neurobiology's early years it was common sense that neuron should communicate with each other remotely, only by means of their electrical activity. In this context, the theory of so-called \emph{ephaptic connection} was forged in the 1940s by A.\ Arvanitaki, Nobel laureate B.\ Katz, and H.\ O.\ Schmitt, cf.~\cite{Arv42,KatSch40}. Such a theory was thought to be surpassed after the newly invented electron microscopes allowed in 1954 to finally prove the existence of chemical synapses.

Although synaptical connections are ultimately stronger and more common, more recent experiments have however found evidence of ephaptic effects in several animals and even in human patients. While experiments have been conducted in real neuronal networks, to the best of our knowledge mathematical models of ephaptic connections have only been treated in \cite{Hol98}. Though, there is some literature for the special case of bundles of (synaptically) unconnected nerve fibers of infinite length, cf.~\cite{HolKoc99,BinEilSco01,BokLaaBli01} and references therein.

Although the derivation in the quoted articles is different, the mathematical models presented in~\cite[\S 4]{HolKoc99}, \cite[\S 4]{BinEilSco01}, and \cite{BokLaaBli01} are comparable. Possibly up to linearization, they describe ephaptic interaction within a myelinated nerve fiber of $m$ axons of infinite lengths by a system of diffusion equations of the form
\begin{equation}\label{problem2}
\left\{
\begin{array}{rcll}
\dot{u}_1(t,x)&=& \sum_{j=1}^m (c_{1j} u'_j(t,\cdot))'(x), \quad & t\geq 0,\; x\in{\mathbb R},\\
&\vdots&\\
\dot{u}_m(t,x)&=& \sum_{j=1}^m (c_{mj} u'_j(t,\cdot))'(x), \quad & t\geq 0,\; x\in{\mathbb R},\\
\end{array}
\right.
\end{equation}
where $u_i(t,x)$ is the electric potential of the $i$\textsuperscript{th} axon at space $x$ and time $t$. The  coefficients $(c_{ij})$ are positive constants that represent the ephaptic effect on the $i$\textsuperscript{th} axon due to the activity of the $j$\textsuperscript{th} one. We emphasize that the mutual interactions and therefore the matrix $(c_{ij})$ are in general non-symmetric.

Whenever potential transmission in neuronal networks is mathematically modelled, neurobiologists usually assume that some form of Kirchhoff law holds in the nodes, as well as continuity of potential. In the easiest linear case, this amounts to saying that in each node the total incoming electric flow equals the total outgoing one, possibly up to some form of dissipation, cf.~\cite{MajEvaJac93}. As we have seen in Remark~\ref{kirchgen}, the natural generalization of Kirchhoff node conditions to the case of strongly coupled network equations is given by~\eqref{kkk}.

This motivates us to consider~\eqref{problem} as a model for transmission of potential in (passive) nerve fibers where ephaptic effects hold. The following results allows to easily discuss also the computationally hard case of numerous contiguous neurons.

\begin{prop}
If the coefficients $c_{ij}$ satsify 
\begin{equation}\label{gersapplic1}
c_{ii}> \sum _{j\not= i} \frac{|c_{ij}+c_{ji}|}{2}, \qquad i=1,\ldots,m,
\end{equation}
then the initial value problem associated with~\eqref{problem2} is well-posed.
\end{prop}

\begin{proof}
By the results of Section~2, the initial value problem is well-posed if the coefficient matrix $C$ is coercive. By Gershgorin's Circle Theorem, we directly obtain that~\eqref{gersapplic1} implies coercivity of the matrix $C$, and the assertion follows by Corollary~\ref{wp}.
\end{proof}
 
The coefficients $(c_{ij})$ are phenomenological constants that have to be determined experimentally. As already observed in~\cite[\S~4.1]{CarMug07}, the model proposed in~\cite{HolKoc99} (i.\,e., $c_{ij}\equiv c$ for all $i,j$) seems to be ill-posed in the light of Remark~\ref{ahiahiahi}, whereas in the models proposed in~\cite{BinEilSco01}--\cite{BokLaaBli01} the possibility to apply Corollary~\ref{wp} depends on the values given to the coupling parameters.

In all models of ephaptic coupling considered above, the coefficients are assumed to satisfy $\sum_{i=1}^m c_{ij}=\mathrm{const}_1$ for all $j$ and $\sum_{j=1}^m c_{ij}=\mathrm{const}_2$ for all $i$. Then by Theorem~\ref{mittelunggraph} one can say that a necessary condition for the subspace of pointwise equal functions to be invariant under the action of $(e^{tA})_{t\geq 0}$ is that the neuronal network is either bipartite or Eulerian. In fact, assuming for the sake of simplicity that no dissipation happens in the nodes (i.\,e., $M=0$), one deduces that there exists two function $C_1,C_2: [0,1] \to \mathbb C$ such that for all $x \in [0,1]$ $\sum_{j=1}^m c_{ij}(x)=C_1(x)$ for all $i$ and $\sum_{i=1}^m c_{ij}(x)=C_2(x)$ for all $j$. 

Observe that by Proposition~\ref{contractive} even if the system is governed by a contractive semigroup in $X^2$ (which is the case if $M$ is dissipative), no contractivity property holds with respect to the norms $\|\cdot\|_1$ and $\|\cdot\|_\infty$ unless $C$ is diagonal. In other words, the system's potential may increase both globally and locally, as soon as ephaptic effects are actually considered.

\subsection{Quantum graphs}\label{quantum}
Consider a finite network of thin waveguides $\me_1,\ldots,\me_m$ of (possibly different) lengths $\ell_1,\ldots,\ell_m$. Discussing the propagation of wave functions, i.\,e., studying the evolution of a system of Schr\"odinger equations
$$i\hslash \frac{\partial v_j}{\partial t}(t,x)=\frac{\partial^2 v_j}{\partial x^2}(t,x),\qquad x\in(0,\ell_j),\; t\in{\mathbb R},$$
over such linear structures -- usually called \emph{quantum graphs} -- has become a relevant topic in recent years,
see e.g.~\cite{KosSch99}--\cite{GutSmi01}--\cite{Kuc02} and references therein. Kirchhoff or more general self-adjoint conditions are usually imposed in the nodes of quantum graphs.

In order to define an Hamiltonian associated with the quantum graph, observe that after a change of coordinates the above equation reads
$$\frac{\partial u_j}{\partial t}(t,x)=\frac{-i}{\hslash\ell^2_j} \frac{\partial^2 u_j}{\partial x^2}(t,x),\qquad x\in(0,1),\; t\in{\mathbb R}.$$
The Hamiltonian is thus given by $iA$, where $(A,D(A))$ is the operator introduced in~\eqref{operator}--\eqref{domain} and associated with the form $a$. Here we are considering coefficients
$$c_{ij}=\left\{\begin{array}{ll}
\frac{1}{\hslash\ell^2_j}\qquad &\hbox{if }i=j,\\
0 & \hbox{otherwise}.
               \end{array}\right.$$
Thus, the operator $A$ is self-adjoint if and only if the ephaptic coupling and nodal coefficient
matrices $C(x)$, $x\in[0,1]$, and $M$ are both self-adjoint, which we assume throughout.
Then by Stones's theorem $iA$ generates a  unitary group that governs the evolution
on the quantum graph.
As in classical field theory, we introduce the action functional $\mathcal S$ for the time evolution of the quantum graph (resp., of the parabolic problem), which is defined as
$${\mathcal S}(\psi)= \int_0^T \sum_{j=1}^m \int_0^1 \left(i \psi_j \overline{\dot{\psi_j}} +\frac{1}{2\hbar \ell^2_j} |\psi_j'|^2\right) dxdt$$
(resp., as
$${\mathcal S}(\psi)= \int_0^T \sum_{j=1}^m \int_0^1 \left(\psi_j \overline{\dot{\psi_j}} +\frac{1}{2\hbar \ell^2_j} |\psi_j'|^2\right) dxdt),$$
i.\,e., 
${\mathcal S}(\psi)=\int_0^T \left(i(\psi|\dot{\psi})_{X^2}+a(\psi,\psi) \right) dt$ 
(resp.,  $\mathcal{S}(\psi) = \int_0^T \left( (\psi|\dot{\psi})_{X^2}+a(\psi,\psi) \right)dt$). 
Here we have implicitely assumed that $\psi \in C^1([0,T]), X^2) \cap C([0,T],V)$ for an arbitrary $T>0$.
Our aim is to discuss symmetry property of the system, in the following sense.

\begin{defi}\label{sym} $ $
\begin{enumerate}
\item We call a $C_0$-group $(U(s))_{s\in\mathbb R}$ on $X^2$ a \emph{symmetry group}
for the system of Schr\"odinger equations (parabolic equations) over the network
if ${\mathcal S}(\psi)={\mathcal S}(\mathcal U(s)\psi)$ for all $s\in\mathbb R$, where
$
(\mathcal U(s) \psi)(t):= U(s)\psi(t)$, $t\in [0,T].$
\item We say that a bounded linear operator $\mathcal P$ on $X^2$ \emph{reflects a symmetry} of the parabolic network problem if
$\mathcal{P}_K e^{tA} = e^{tA} \mathcal{P}_K$ for all $t \ge 0$, i.\,e., if projecting the initial value and then studying the corresponding time evolution is equivalent to projecting the solution curve of the original problem.
\end{enumerate}
\end{defi}

Since $\mathcal U(s)$ does not act on the time variable, one sees that due to the time continuity of $\psi$ a self-adjoint bounded linear operator on $X^2$, i.\,e., an \emph{observable} $\mathcal P$ of the physical system, is the infinitesimal generator of a symmetry group $(e^{is\mathcal P})_{s\in\mathbb R}$ if and only if it satisfies $a(\psi,\psi)=a(e^{is\mathcal P}\psi,e^{is\mathcal P}\psi)$ for all $s\in\mathbb R$ and all $\psi \in V$.
	

We consider the case of  a closed linear subspace $\mathcal Y$ constructed as in~\eqref{projy} and discuss observables ${\mathcal P}={\mathcal P}_K$ given by orthogonal projections of the state space $X^2:=(L^2(0,1))^m$ onto $\mathcal Y$ that satisfy~\eqref{projdef}. A justification for the use of the term ``symmetry'' in Definition~\ref{sym} is given in the following.

\begin{prop}\label{gauge}
Let $\mathcal{P}_K$ be an observable of the system as defined above, and assume $K$ to be admissible. The following assertions are equivalent.
\begin{enumerate}[(a)]
\item\label{reflsymm}
	The projection $\mathcal{P}_K$ reflects a symmetry of the network parabolic problem.
\item\label{invariance}
	The subspace $\mathcal{Y} = \Range \mathcal{P}_K$ is invariant under $(e^{tA})_{t \ge 0}$.
\item\label{formproj}
	$a(\mathcal{P}_K \psi, \psi) = a(\mathcal{P}_K\psi, \mathcal{P}_K\psi)$ for all $\psi \in V$.
\item\label{formgroup}
	The projection $\mathcal{P}_K$ generates a symmetry group of the parabolic network equation,
	i.\,e., $a(\psi, \psi) = a(e^{is\mathcal{P}_K}\psi, e^{is\mathcal{P}_K}\psi)$ for all $s \in \mathbb{R}$ and all $\psi \in V$.
\item\label{formgroupi}
	The projection $\mathcal{P}_K$ generates a symmetry group of the network Schr\"odinger equation.
\item\label{invariancei}
	The subspace $\mathcal{Y} = \Range \mathcal{P}_K$ is invariant under $(e^{itA})_{t \in \mathbb{R}}$.
\end{enumerate}
\end{prop}

\begin{proof}
Note that the invariance of $\mathcal{Y} = \Range \mathcal{P}_K$ under the action
of $(e^{tA})_{t\geq 0}$ is equivalent to
\[\tag{\ref{invariance}'}\label{invariance2}
	\mathcal{P}_K e^{tA} \mathcal{P}_K = e^{tA} \mathcal{P}_K \text{ for all } t\geq 0.
\]

``\eqref{reflsymm} $\Rightarrow$ \eqref{invariance2}''
This is obvious, since $\mathcal{P}_K^2 = \mathcal{P}_K$.

``\eqref{invariance2} $\Rightarrow$ \eqref{reflsymm}''
Since $\mathcal{P}_K$ and $e^{tA}$ are self-adjoint,
\[
\mathcal{P}_K e^{tA} = \left( e^{tA} \mathcal{P}_K \right)^\ast
= \left( \mathcal{P}_K e^{tA} \mathcal{P}_K \right)^\ast
= \mathcal{P}_K e^{tA} \mathcal{P}_K
= e^{tA} \mathcal{P}_K.
\]

``\eqref{invariance} $\Leftrightarrow$ \eqref{formproj}''
By Theorem~\ref{oulin}, \eqref{invariance} is equivalent to $a(\mathcal{P}_K f, (\Id - \mathcal{P}_K) f) = 0$
for every $f \in X^2$. But this is \eqref{formproj}.

``\eqref{formproj} $\Leftrightarrow$ \eqref{formgroup}''
Since $\mathcal{P}_K$ is a projection,
\[
	e^{z\mathcal{P}_K}=\sum_{j=0}^\infty\frac{z^j}{j!}{\mathcal P}_K^j=\sum_{j=1}^\infty \frac{z^j}{j!}{\mathcal P}_K+\Id=(e^{z}-1){\mathcal P}_K+\Id.
\]
Using this representation we see that
\begin{eqnarray*}
a\left(e^{is\mathcal{P}_K}\psi, e^{is\mathcal{P}_K}\psi\right)
& = & a\left( (e^{is} - 1) \mathcal{P}_K \psi, (e^{is} - 1) \mathcal{P}_K \psi \right)\\
&&+ 2 \Real a\left( (e^{is} - 1) \mathcal{P}_K \psi, \psi \right)
+ a(\psi, \psi) \\
& = &
a\left( \mathcal{P}_K \psi, \mathcal{P}_K \psi \right)
- 2\Real e^{is} a\left( \mathcal{P}_K \psi, \mathcal{P}_K \psi \right)
+ a\left( \mathcal{P}_K \psi, \mathcal{P}_K \psi \right) \\ & &
+ 2 \Real (e^{is} - 1) a\left( \mathcal{P}_K \psi, \psi \right)
+ a(\psi, \psi) \\
& = & 2 a\left( \mathcal{P}_K \psi, \mathcal{P}_K \psi \right) - 2 a\left( \mathcal{P}_K \psi, \psi \right) + a(\psi, \psi).
\end{eqnarray*}
Thus \eqref{formgroup} is equivalent to
$2 a\left( \mathcal{P}_K \psi, \mathcal{P}_K \psi \right) - 2 a\left( \mathcal{P}_K \psi, \psi \right) = 0$
{ for every } $\psi \in V$, which is \eqref{formproj}.

``\eqref{formgroup} $\Leftrightarrow$ \eqref{formgroupi}''
Both statements are equivalent to $a(\psi, \psi) = a(e^{is\mathcal{P}_K}\psi, e^{is\mathcal{P}_K}\psi)$ for all $s \in \mathbb{R}$
and all $\psi \in V$, since $e^{is\mathcal{P}_K}$ is an unitary operator that commutes
with the time derivative.

``\eqref{invariance} $\Rightarrow$ \eqref{invariancei}''
After rescaling we may assume that $(e^{tA})_{t\geq 0}$ is contractive.
It is known that the invariance of ${\mathcal Y}$ under $(e^{tA})_{t\geq 0}$ is equivalent to invariance of $\mathcal Y$ under $R(\lambda,A)$ for all $\lambda\in\mathbb R$ large enough, see e.g.~\cite[Prop.~2.1]{Ouh04}. On the other hand, $(e^{itA})_{t\in\mathbb R}$ is also a (unitary, hence contractive) $C_0$-(semi)group that satisfies $e^{itA}{\mathcal Y}\subset {\mathcal Y}$ if and only if $R(\lambda,iA){\mathcal Y}\subset {\mathcal Y}$ for $\lambda$ large enough, i.\,e., if and only if  $ iR\left(\frac{\lambda}{i},A\right){\mathcal Y}\subset {\mathcal Y}$ for $\lambda$ large enough.
In fact, the resolvent set of $A$ contains an open sector of $\mathbb C$ which  contains $\overline{\Sigma}$, i.\,e., it contains the closed right half plane (with the possible exception of the origin). Then, for any $\lambda_0,\mu\in \Sigma$ such that $|\mu-\lambda_0| \leq \|R(\lambda_0,A)\|^{-1}$ it is possible to develop the resolvent operator $R(\mu,A)$ as a power series centered at $\lambda_0$, i.\,e., 
$$
R(\mu, A)= \sum\limits_{n=0}^\infty (\lambda_0-\mu)^n R(\lambda_0,A)^{n+1}.
$$
Let now ${\mathcal Y}$ be invariant under $(e^{tA})_{t\geq 0}$. Then ${\mathcal Y}$ is invariant under $R(\lambda_0,A)$ for some $\lambda_0$, i.\,e., $R(\lambda_0,A)y\in {\mathcal Y}$ for all $y\in {\mathcal Y}$. Since ${\mathcal Y}$ is a closed linear subspace, one obtains 
$$R(\mu,A)y=\sum_{n=0}^\infty (\lambda_0-\mu)^n R(\lambda_0,A)^{n+1}y \in {\mathcal Y}\qquad\hbox{for all } y\in {\mathcal Y} \hbox{ and }|\mu-\lambda_0| \leq \|R(\lambda_0,A)\|^{-1},$$
and therefore ${\mathcal Y}$ is invariant under $R(\mu,A)$. This shows that the subset of the resolvent set for which $R(\lambda,A)\mathcal Y \subset Y$ is open. Moreover, it is relatively closed since $\mathcal Y$ is closed. As a consequence, $\mathcal Y$ is invariant under $R(\mu,A)$ for all $\mu$ in the unbounded connected component of the spectrum containing $\lambda_0$, and therefore also for all $i \lambda$, $\lambda \in \mathbb R$ large enough. By the representation of the semigroup in terms of the resolvent this shows that ${\mathcal Y}$ is invariant under the unitary group $(e^{itA})_{t\in\mathbb R}$.

``\eqref{invariancei} $\Rightarrow$ \eqref{invariance}''
This can be proved in the same spirit as the implication ``\eqref{invariance} $\Rightarrow$ \eqref{invariancei}''.
\end{proof}

\begin{rem}
Careful examination of the proof above shows that $a(f, f) = a\left(e^{it\mathcal{P}_K}f, e^{it\mathcal{P}_K}f\right)$
and admissibility of $K$ are equivalent to 
\begin{equation}\label{synnonsa}
\mathcal{P}_K e^{tA} = e^{tA^\ast} \mathcal{P}_K\qquad \hbox{for all }t\geq 0,
\end{equation}
even if $A$ is not self-adjoint. Definition~\ref{sym} can thus be generalized by saying that \emph{${\mathcal P}_K$ reflects a symmetry of a (possibly non-self-adjoint) parabolic network problem} if~\eqref{synnonsa} holds.
\end{rem}

Thus, $\mathcal Y$ is invariant under the action of $(e^{itA})_{t\geq 0}$ if and only if the associated orthogonal
projection $\mathcal P$ is admissible and the orthogonality condition is satisfied by $M$, i.\,e., if and only
if $\mathcal Y$ is invariant under the action of $(e^{tA})_{t \ge 0}$.
In particular, for a star graph $\mG$ Proposition \ref{star}.(1) yields that there are nontrivial invariant subspaces
of the above form if and only if there is a pair of edges with the same length.

\bibliographystyle{plain} 
\bibliography{/home/stefano/Desktop/matematica/dissertation/literatur.bib}
\end{document}